\documentclass{article}

\DeclareMathAlphabet{\mathpzc}{OT1}{pzc}{m}{it}

\usepackage{graphicx}
\usepackage{amssymb}
\usepackage{amsmath}
\usepackage{amsthm}
\usepackage{bm}

\usepackage[initials]{amsrefs}

\numberwithin{equation}{section}
\allowdisplaybreaks

\newtheorem{theorem}{Theorem}
\newtheorem{remark}{Remark}

\date{}

\author{Enzo Orsingher, Federico Polito, \& Ludmila Sakhno}

\title{Fractional Non-Linear, Linear and Sublinear Death Processes}

\begin{document}
	\maketitle
	\begin{abstract}
		This paper is devoted to the study of a fractional version of non-linear
		$\mathpzc{M}^\nu(t)$, $t>0$, linear $M^\nu (t)$, $t>0$ and sublinear
		$\mathfrak{M}^\nu (t)$, $t>0$ death processes. Fractionality is introduced by replacing the usual
		integer-order derivative in the difference-differential equations governing
		the state probabilities, with the fractional derivative understood in the sense of
		Dzhrbashyan--Caputo.
		We derive explicitly the state probabilities of the three death processes and
		examine the related probability generating functions and mean values. A
		useful subordination relation is also proved, allowing us to express the death processes as
		compositions of their classical counterparts with the random time process
		$T_{2 \nu} (t)$, $t>0$. This random time has one-dimensional distribution
		which is the folded solution to a Cauchy problem of the fractional diffusion
		equation.
	\end{abstract}

	\section{Introduction}

		We assume that we have a population of $n_0$ individuals or objects. The components of
		this population might be the set of healthy people during an epidemic or the
		set of items being sold in a store, or even, say, melting ice pack blocks.
		However even a coalescence of particles can be treated in this same manner,
		leading to a large ensemble of physical analogues suited to the method.
		The main interest is to model the fading process of these objects
		and, in particular, to analyse how the size of the population decreases.

		The classical death process is a model describing this type of phenomena and, its linear
		version is analysed in \ocite{bailey}, page 90. The most interesting feature of the
		extinguishing population is the probability distribution
		\begin{equation}
			\label{probab}
			p_k (t) = \text{Pr} \left\{ M(t) = k \mid M(0) = n_0 \right\}, \qquad t>0,
			\: 0 \leq k \leq n_0,
		\end{equation}
		where $M(t)$, $t>0$ is the point process representing the size of the population at time
		$t$.
		If the death rates are proportional to the population size, the process is called
		\emph{linear} and the probabilities \eqref{probab} are solutions to the initial-value
		problem
		\begin{equation}
			\label{cauchy}
			\begin{cases}
				\frac{d}{dt} p_k (t) = \mu (k+1) p_{k+1} (t) - \mu k p_k (t),
				& 0 \leq k \leq n_0, \\
				p_k(0) =
				\begin{cases}
					1, & k=n_0, \\
					0, & 0 \leq k < n_0,
				\end{cases}
			\end{cases}
		\end{equation}
		with $p_{n_0+1}(t) = 0$.

		The distribution satisfying \eqref{cauchy} is
		\begin{equation}
			\label{class}
			p_k (t) = \binom{n_0}{k} e^{-\mu k t} \left( 1-e^{-\mu t} \right)^{n_0-k},
			\qquad 0 \leq k \leq n_0.
		\end{equation}
		The equations \eqref{cauchy} are based on the fact that the death rate of each component of the
		population is proportional to the number of existing individuals.

		In the non-linear case, where the death rates are $\mu_k$, $0 \leq k \leq n_0$,
		equations \eqref{cauchy} must be replaced by
		\begin{equation}
			\label{cauchy-nl}
			\begin{cases}
				\frac{d}{dt} \mathpzc{p}_k (t) = \mu_{k+1} \mathpzc{p}_{k+1} (t) - \mu_k
				\mathpzc{p}_k (t),
				& 0 \leq k \leq n_0, \\
				\mathpzc{p}_k(0) =
				\begin{cases}
					1, & k=n_0, \\
					0, & 0 \leq k < n_0.
				\end{cases}
			\end{cases}
		\end{equation}

		In this paper we consider fractional versions of the processes
		described above, where fractionality is obtained by substitution of the integer-order
		derivatives appearing in \eqref{cauchy} and \eqref{cauchy-nl},
		with the fractional derivative called Caputo or
		Dzhrbashyan--Caputo derivative,
		defined as follows
		\begin{equation}
			\label{caputo}
			\begin{cases}
				\frac{d^\nu f \left( t \right) }{dt^\nu} = \frac{1}{\Gamma \left( 1- \nu \right)}
				\displaystyle \int_0^t \frac{ f'
				\left( s \right)}{\left( t-s \right)^\nu} \, ds,
				& 0 < \nu < 1, \\
				f' \left( t \right), & \nu = 1.
			\end{cases}
		\end{equation}
		The main advantage of the Dzhrbashyan--Caputo fractional derivative over the usual
		Riemann--Liouville fractional derivatives is that the former requires only integer-order
		derivatives in the initial conditions.

		The fractional derivative operator is vastly present in the physical and mathematical
		literature. It appears for example in generalisations of diffusion-type differential
		equations (see \ocite{wyss}, \ocite{wyss2}, \ocite{nigmatullin} and \ocite{mainardi}),
		hyperbolic equations such as telegraph equation (see \ocite{orsbeg}),
		reaction-diffusion equations (see \ocite{saxena1}), or in the study of
		continuous time random walks (CTRW) scaling limits (see \ocite{kol1}, \ocite{meer}).
		Fractional calculus has also been considered by some authors to describe
		cahotic Hamiltonian dynamics in low dimensional systems (see e.g.\ \ocite{zas},
		\ocite{zas2}, \ocite{saxena2}, \ocite{saxena3} and \ocite{saxena4}).
		For a complete review of fractional kinetics the reader can consult \ocite{zas3}
		or the book by Zaslavsky \ocite{zas4}.
		In the literature are also present fractional generalisations of point processes,
		such as the Poisson process (see \ocite{repin}, \ocite{laskin}, \ocite{scalas},
		\ocite{cahoy}, \ocite{sibatov} and \ocite{orsbeg2}) and the birth and birth-death processes
		(see \ocite{cahoy2}, \ocite{pol}, \ocite{pol2}). Fractional models are also used in
		other fields, for example finance (\ocite{scalas2}, \ocite{scalas3}).

		The population size is governed by
		\begin{equation}
			\label{cauchy-fnl}
			\begin{cases}
				\frac{d^\nu}{dt^\nu} \mathpzc{p}_k (t) = \mu_{k+1} \mathpzc{p}_{k+1} (t) - \mu_k
				\mathpzc{p}_k (t),
				& 0 \leq k \leq n_0, \\
				\mathpzc{p}_k^\nu(0) =
				\begin{cases}
					1, & k=n_0, \\
					0, & 0 \leq k < n_0,
				\end{cases}
			\end{cases}
		\end{equation}
		and is denoted by $\mathpzc{M}^\nu (t)$, $t>0$.

		Let us assume that a crack has the form of a process $T_{2\nu}(t)$, $t>0$. For $\nu=1/2$, this
		coincides with a reflecting Brownian motion and has been described and derived in
		\ocite{kunin}. For $\nu \neq 1/2$, the process $T_{2 \nu} (t)$, $t>0$,
		can be identified with a stable process (see for details on this point \ocite{orsbeg3}).
		The ensemble of $n_0$ particles moves on the fracture and, at the same time,
		undergoes a decaying process which respects the same probabilistic rules of the usual
		death process.
		For the number of existing particles, we have therefore
		\begin{equation}
			\label{uno}
			\mathpzc{p}^{\nu}_{k}(t) = \int_0^\infty \mathpzc{p}_{k}(s) \text{Pr} \left\{T_{2\nu}(t)
			\in ds \right\}.
		\end{equation}
		We observe that
		\begin{equation}
			\text{Pr} \left\{ T_{2\nu}(t) \in ds \right\} = q(s, t) ds,
		\end{equation}
		is a solution to
		\begin{equation}
			\frac{\partial^{2\nu} }{\partial t^{2\nu}} q(s,t)
			= \frac{\partial^2 }{\partial s^2} q(s,t), \qquad s>0,\:t>0,
		\end{equation}
		with the necessary initial conditions.
		Furthermore we recall that
		\begin{equation}
			\int_0^\infty e^{-z t} q(s, t) dt
			= z^{\nu-1} e^{-z^{\nu} s}, \qquad z >0, \: s > 0.
		\end{equation}
		The distribution $q(s,t)$ is also a solution to
		\begin{equation}
			\label{due}
			\frac{\partial^{\nu}}{\partial t^{\nu}} q(s,t)
			= -\frac{\partial}{\partial s} q(s,t), \qquad s>0,
		\end{equation}
		as can be ascertained directly. If we take the fractional derivative in \eqref{uno} we get
		\begin{align}
			\frac{d^{\nu}}{d t^{\nu}} \mathpzc{p}^{\nu}_k (t)
			& = \int_0^\infty \mathpzc{p}_{k}(s) \frac{\partial^{\nu}}{\partial t^{\nu}}
			\text{Pr} \left\{ T_{2\nu}(t) \in ds \right\} \\
			& = - \int_0^\infty \mathpzc{p}_{k}(s) \frac{\partial}{\partial s}
			q(s,t) ds \notag \\
			& = - q(s,t) \mathpzc{p}_{k}(s) \big|_{0}^{\infty}
			+ \int_0^\infty \frac{d \mathpzc{p}_{k}(s)}{d s}
			q(s,t) ds \notag \\
			& = \int_0^\infty \left[-\mu_k \mathpzc{p}_k(s) + \mu_{k+1} \mathpzc{p}_{k+1}(s)
			\right] q(s,t) ds \notag \\
			& = - \mu_k \mathpzc{p}_k^{\nu}(t) + \mu_{k+1} \mathpzc{p}^{\nu}_{k+1}(t). \notag
		\end{align}
		This shows that replacing the time derivative with the fractional derivative corresponds
		to considering a death process (annihilating process) on particles displacing on a crack.

		We now give some details about \eqref{due}.
		By taking the Laplace transform of both members of \eqref{due} we have that
		\begin{align}
			\int_0^\infty e^{-z t} \frac{\partial^{\nu}}{\partial t^{\nu}} q(s,t)
			dt & = - \frac{\partial}{\partial s} \left( \int_0^\infty e^{-z t}
			q(s,t) dt \right) \\
			& = - \frac{\partial}{\partial s} \left( z^{\nu -1} e^{-s z^{\nu}} \right) =
			z^{2\nu -1} e^{-s z^{\nu}}. \notag
		\end{align}
		Furthermore,
		\begin{align}
			\int_0^\infty e^{-z t} \frac{\partial^{\nu}}{\partial t^{\nu}} q(s,t) dt
			& = z^{\nu} \int_0^\infty e^{-z t} q(s,t) dt - z^{\nu -1} q(s,0) \\
			& = z^{\nu} \left( z^{\nu -1} e^{- s z^{\nu}}\right) - z^{\nu -1 } \delta(s), \notag
		\end{align}
		and therefore for $s>0$ this establishes that $q(s,t)$ solves equation \eqref{due}.
		We note that a gas particle moving on a fracture has inspired to different authors
		the iterated Brownian motion (see \ocite{deblassie}).

		The distribution
		\begin{equation}
			\mathpzc{p}_k^\nu (t) =
			\text{Pr} \left\{ \mathpzc{M}^\nu (t) =
			k \mid \mathpzc{M}^\nu (0) = n_0 \right\}, \qquad 0 \leq k \leq n_0,
		\end{equation}
		is obtained explicitly
		and reads
		\begin{equation}
			\label{dist-nl}
			\mathpzc{p}_k^\nu (t) =
			\begin{cases}
				E_{\nu,1} (-\mu_{n_0} t^\nu), & k = n_0, \\
				\textstyle\prod\limits_{j=k+1}^{n_0} \mu_j
				\displaystyle\sum_{m=k}^{n_0} \frac{ E_{\nu,1} (-
				\mu_m t^\nu)}{\textstyle\prod\limits_{\substack{
				h=k \\ h \neq m}}^{n_0} \left( \mu_h -
				\mu_m \right) }, & 0 < k < n_0, \\
				1 - \displaystyle\sum_{m=1}^{n_0} \displaystyle\prod_{\substack{
				h=1 \\ h \neq m}}^{n_0} \left( \frac{\mu_h}{\mu_h-\mu_m} \right)
				E_{\nu,1} (- \mu_m t^\nu ), & k=0, \: n_0 > 1.
			\end{cases}
		\end{equation}
		Obviously, for $k=0$, $n_0=1$,
		\begin{equation}
			\mathpzc{p}_0^\nu (t) = 1-E_{\nu,1} (- \mu_1 t^\nu).
		\end{equation}

		The Mittag-Leffler functions appearing in \eqref{dist-nl} are defined as
		\begin{equation}
			\label{mittag}
			E_{\nu,\gamma} \left( x \right) = \sum_{h=0}^\infty
			\frac{x^h}{\Gamma \left( \nu h+\gamma \right)}, \qquad
			x \in \mathbb{R}, \quad \nu, \gamma > 0 .
		\end{equation}
		For $\nu=\gamma=1$, $E_{1,1} (x) = e^x$ and formulae \eqref{dist-nl} provide the explicit
		distribution of the classical non-linear death process.

		For $\mu_k = k \mu$ the distribution of the fractional linear death process can be
		obtained either directly by solving the Cauchy problem \eqref{cauchy-fnl}
		with $\mu_k=k \cdot \mu$ and $p_{n_0+1}(t)=0$, or
		by specialising \eqref{dist-nl} resulting in the following form
		\begin{equation}
			\label{distr}
			p_k^\nu (t) = \binom{n_0}{k} \sum_{r=0}^{n_0-k} \binom{n_0-k}{r} (-1)^r E_{\nu,1}
			(- (k+r) \mu t^\nu ).
		\end{equation}

		A technical tool necessary for our manipulations is the Laplace transform of
		Mittag-Leffler functions which we write here for the sake of completeness:
		\begin{equation}
			\int_0^\infty e^{-zt} t^{\gamma - 1} E_{\nu,\gamma} (\pm \vartheta t^\nu) dt
			= \frac{z^{\nu-\gamma}}{z^\nu \mp \vartheta}, \qquad \mathfrak{R}(z) >
			|\vartheta|^\frac{1}{\nu}.
		\end{equation}

		Another special case is the so-called fractional sublinear death process
		(for sublinear birth processes consult \ocite{donnelly}) where the death rates
		have the form $\mu_k = \mu (n_0+1-k)$. In the sublinear process, the annihilation
		of particles or individuals accelerates with decreasing population size.

		The distribution $\mathfrak{p}_k^\nu (t)$, $0 \leq k \leq n_0$ of the fractional
		sublinear death process $\mathfrak{M}^\nu (t)$, $t>0$, is strictly related
		to that of the fractional linear birth process $N^\nu (t)$, $t>0$
		(see, for details on this point, \ocite{pol}):
		\begin{equation}
			\text{Pr} \left\{ \mathfrak{M}^\nu (t) = 0 \mid
			\mathfrak{M}^\nu (0) = n_0 \right\} = \text{Pr} \left\{ N^\nu (t) > n_0
			\mid N^\nu (0) = 1 \right\}.
		\end{equation}
		In general, the connection between the fractional sublinear death process
		and the fractional linear birth process is expressed by the relation
		\begin{align}
			& \text{Pr} \left\{ \mathfrak{M}^\nu (t) = n_0 - (k-1) \mid
			\mathfrak{M}^\nu (0) = n_0 \right\} \\
			& = \text{Pr} \left\{ N^\nu (t) = k \mid N^\nu (0) = 1 \right\},
			\qquad 1 \leq k \leq n_0. \notag
		\end{align}
		This shows a sort of symmetry in the evolution of fractional linear
		birth and fractional sublinear death processes.

		For all fractional processes considered in this paper, a subordination relationship
		holds. In particular, for the fractional linear death process we can write that
		\begin{equation}
			\label{subsub}
			M^\nu (t) = M (T_{2 \nu} (t) ), \qquad 0 < \nu < 1, \: t>0,
		\end{equation}
		where $T_{2 \nu} (t)$ is a process for which
		\begin{equation}
			\text{Pr} \left\{ T_{2 \nu} (t) \in ds \right\} = q(s,t) ds,
		\end{equation}
		is a solution to the following Cauchy problem (see \ocite{ob})
		\begin{equation}
			\begin{cases}
				\frac{\partial^{2 \nu}}{\partial t^{2 \nu}} q(s,t) =
				\frac{\partial^2}{\partial s^2} q(s,t), & t>0, \: s>0, \\
				\left. \frac{\partial}{\partial t} q(s,t) \right|_{s=0} = 0, \\
				q(s,0) = \delta (s), & 0 < \nu \leq 1,
			\end{cases}
		\end{equation}
		with the additional initial condition
		\begin{equation}
			q_t (s,0) = 0, \qquad 1/2 < \nu \leq 1.
		\end{equation}
		In equation \eqref{subsub}, $M(t)$, $t>0$, represents the classical linear death process.
		Subordination relations of this type are extensively treated in \ocite{orsbeg3}
		and \ocite{koloko}.

		We also show that all the fractional death processes considered below can be viewed
		as classical death processes with rate $\mu \cdot \Xi$, where $\Xi$ is a Wright-distributed
		random variable.

	\section{The fractional linear death process and its properties}

		In this section we derive the distribution of the fractional linear death process
		as well as some interesting related properties and interpretations.
		\begin{theorem}
			The distribution of the fractional linear death process $M^\nu(t)$,
			$t>0$ with $n_0$ initial individuals and death rates $\mu_k = \mu \cdot k$, is given by
			\begin{align}
				\label{distr-linear}
				p_k^\nu (t) & = \text{Pr} \left\{ M^\nu (t) = k \mid M^\nu (0) = n_0 \right\} \\
				& = \binom{n_0}{k} \sum_{r=0}^{n_0-k} \binom{n_0-k}{r} (-1)^r E_{\nu,1}
				(- (k+r) \mu t^\nu ), \notag
			\end{align}
			where $0 \leq k \leq n_0$, $t>0$ and $\nu \in (0,1]$. The function $E_{\nu,1}(x)$ is
			the Mittag-Leffler function previously defined in \eqref{mittag}.

			\begin{proof}
				The state probability $p_{n_0}^\nu (t)$, $t>0$ is readily obtained
				by applying the Laplace transform to equation \eqref{cauchy-fnl},
				with $\mu_k=\mu \cdot k$, and then
				transforming back the result, thus yielding
				\begin{figure}
					\centering
					\includegraphics[scale=0.8,angle=0]{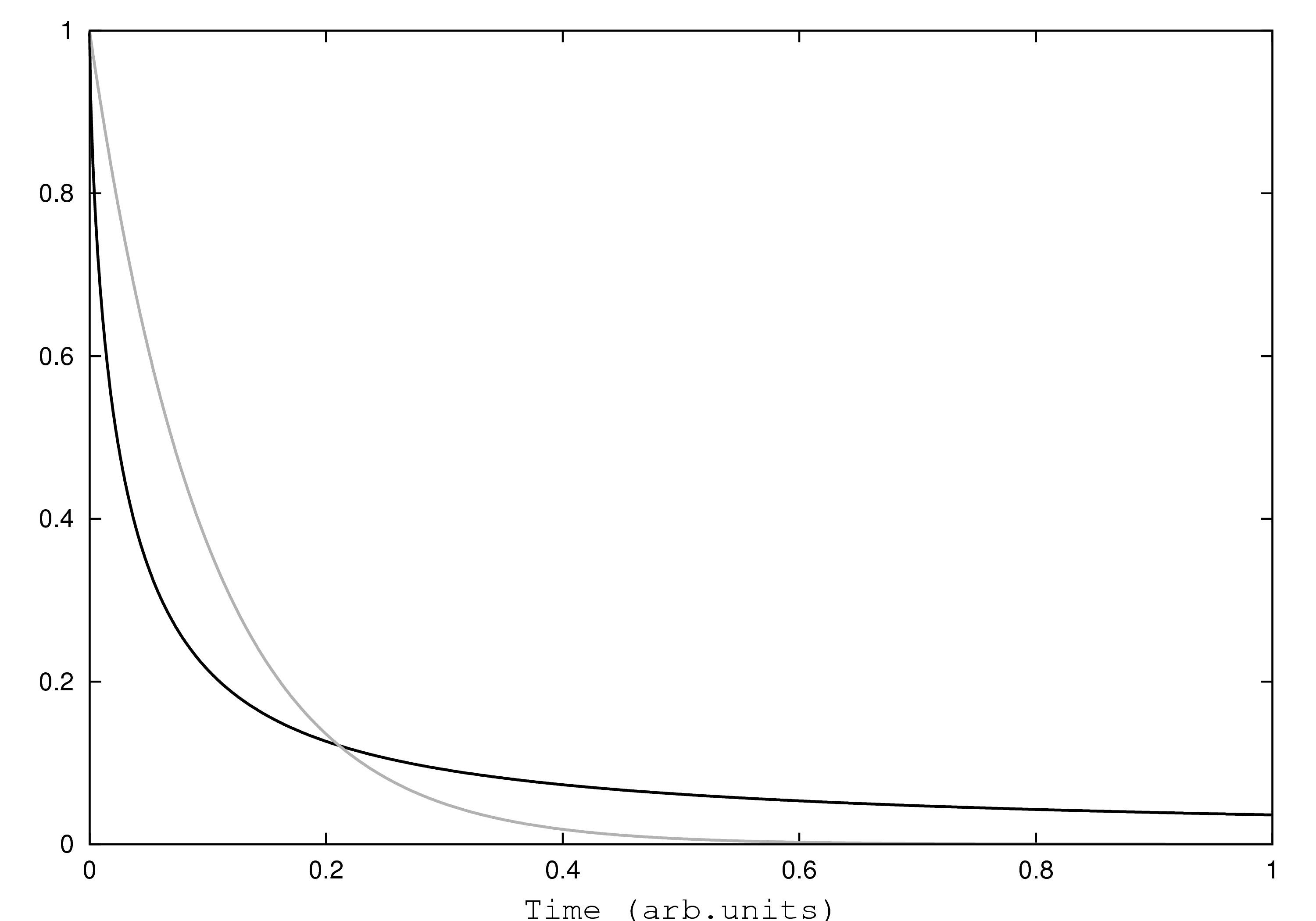}
					\caption{\label{nzero}Plot of $p_{n_0}^{0.7} (t)$ (in black) and
					$p_{n_0}^1 (t)$ (in grey), both with $n_0 = 10$.}
				\end{figure}
				\begin{equation}
					p_{n_0}^\nu (t) = E_{\nu,1} (- n_0 \mu t^\nu), \qquad t>0, \:
					\nu \in (0,1].
				\end{equation}

				When $k=n_0-1$, in order to solve the related differential equation,
				we can write
				\begin{align}
					\label{lap-sec}
					& z^\nu \mathcal{L} \left\{p_{n_0-1}\right\}(z) =
					\mu n_0 \frac{z^{\nu-1}}{z^\nu+n_0 \mu} - \mu (n_0-1) \mathcal{L}
					\left\{p_{n_0-1} \right\} (z) \\
					& \Leftrightarrow \quad \mathcal{L} \left\{ p_{n_0-1} \right\}(z) =
					\mu n_0 z^{\nu -1} \frac{1}{z^\nu +n_0 \mu} \cdot \frac{1}{
					z^\nu + (n_0-1) \mu} \notag \\
					& \Leftrightarrow \quad \mathcal{L} \left\{ p_{n_0-1} \right\}(z) =
					n_0 z^{\nu-1} \left( \frac{1}{z^\nu +(n_0+1) \mu} -
					\frac{1}{z^\nu + n_0 \mu} \right). \notag
				\end{align}

				\pagebreak

				By inverting equation \eqref{lap-sec}, we readily obtain that
				\begin{equation}
					p_{n_0-1}^\nu (t) = n_0 \left( E_{\nu,1} (-(n_0-1)\mu t^\nu) -
					E_{\nu,1} (-n_0 \mu t^\nu ) \right).
				\end{equation}
				\begin{figure}
					\centering
					\includegraphics[scale=0.8,angle=0]{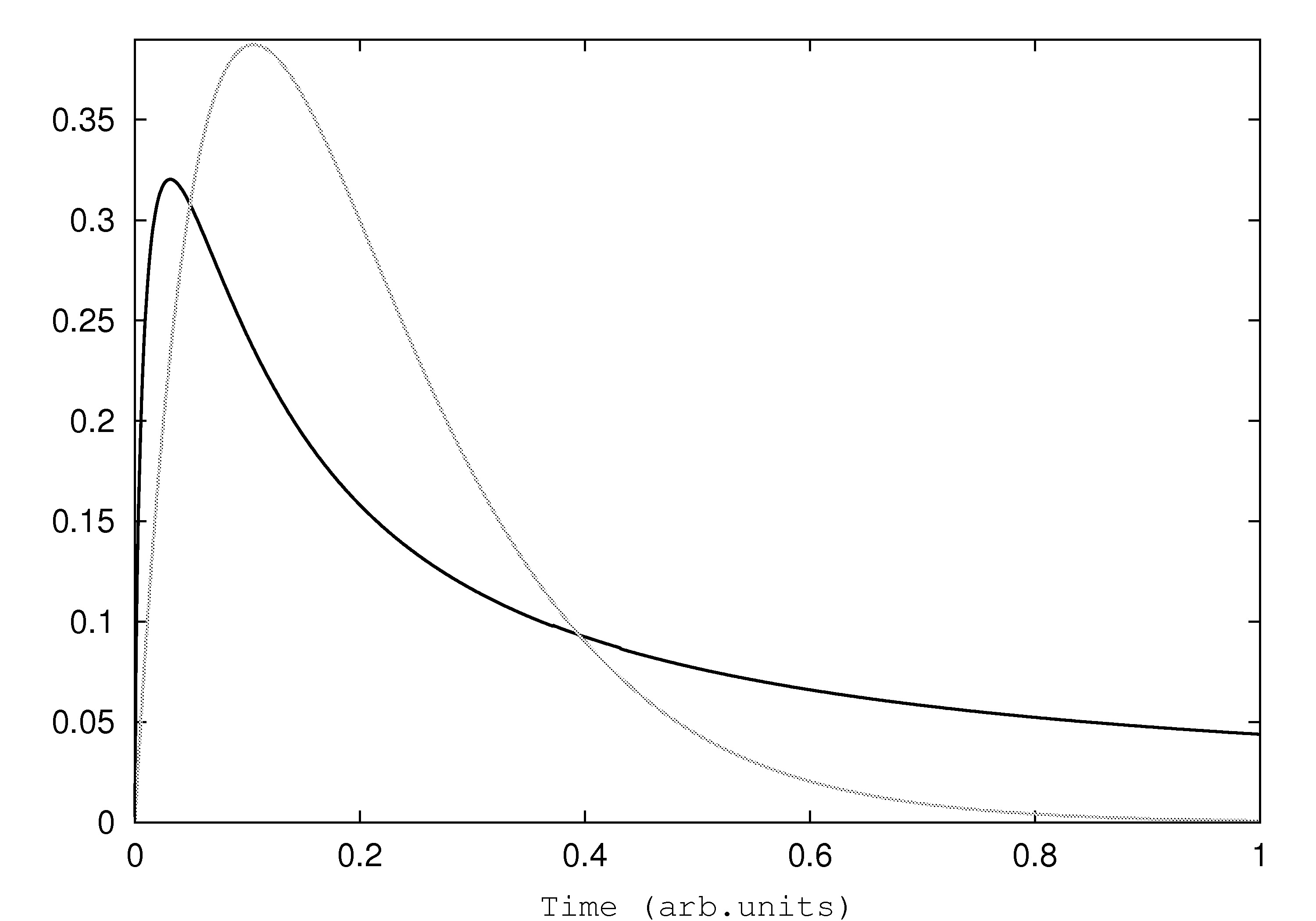}
					\caption{\label{nzeromu}Plot of $p_{n_0-1}^{0.7} (t)$ (in black) and
					$p_{n_0-1}^1 (t)$ (in grey). Here $n_0 = 10$.}
				\end{figure}

				For general values of $k$, with $0 \leq k < n_0$, we must solve the following
				Cauchy problem:
				\begin{align}
					\frac{d^\nu}{dt^\nu} p_k (t) & = \mu (k+1)
					\binom{n_0}{k+1} \\
					& \times \sum_{r=0}^{n_0-k-1} \binom{n_0-k-1}{r}
					(-1)^r E_{\nu,1} (-(k+1+r) \mu t^\nu) - \mu k p_k (t), \notag
				\end{align}
				subject to the initial condition $p_k(0) = 0$ and with $\nu \in (0,1]$.
				The solution can be found by resorting to the Laplace transform, as we see
				in the following.
				\begin{align}
					z^\nu \mathcal{L} \left\{ p_k \right\} (z) & =
					\mu (k+1) \binom{n_0}{k+1} \\
					& \times \sum_{r=0}^{n_0-k-1}
					\binom{n_0-k-1}{r} (-1)^r
					\frac{z^{\nu-1}}{z^\nu+(k+1+r)\mu}
					-\mu k \mathcal{L} \left\{ p_k \right\} (z). \notag
				\end{align}
				The Laplace transform $\mathcal{L} \left\{ p_k \right\} (z)$ can thus
				be written as
				\begin{align}
					\label{lapk}
					& \mathcal{L} \left\{ p_k \right\} (z) \\
					& = \mu (k+1) \binom{n_0}{k+1}
					\sum_{r=0}^{n_0-k-1} \binom{n_0-k-1}{r} (-1)^r
					\frac{z^{\nu-1}}{z^\nu + (k+1+r)\mu} \cdot
					\frac{1}{z^\nu+k\mu} \notag \\
					& = \binom{n_0}{k} \sum_{r=0}^{n_0-k-1} \binom{n_0-k}{r+1}
					(-1)^r z^{\nu-1}
					\left( \frac{1}{z^\nu+k\mu} - \frac{1}{z^\nu+
					(k+1+r)\mu} \right) \notag \\
					& = \binom{n_0}{k} \sum_{j=1}^{n_0-k}
					\binom{n_0-k}{j} (-1)^{j-1} z^{\nu-1} \left( \frac{1}{z^\nu+k\mu}
					- \frac{1}{z^\nu+(k+j)\mu} \right) \notag \\
					& = \binom{n_0}{k} \sum_{j=1}^{n_0-k} \binom{n_0-k}{j}
					(-1)^j \frac{z^{\nu-1}}{
					z^\nu+(k+j)\mu} \notag \\
					& \qquad - \binom{n_0}{k} \frac{z^{\nu-1}}{z^\nu+k\mu}
					\sum_{j=1}^{n_0-k} \binom{n_0-k}{j} (-1)^j \notag \\
					& = \binom{n_0}{k} \sum_{j=1}^{n_0-k} \binom{n_0-k}{j}
					(-1)^j \frac{z^{\nu-1}}{z^\nu+(k+j)\mu} + \binom{n_0}{k}
					\frac{z^{\nu-1}}{z^\nu+k\mu} \notag \\
					& = \binom{n_0}{k} \sum_{j=0}^{n_0-k} \binom{n_0-k}{j}
					(-1)^j \frac{z^{\nu-1}}{z^\nu+(k+j)\mu}. \notag
				\end{align}
			\end{proof}
			By taking now the inverse Laplace transform of \eqref{lapk}, we obtain the
			claimed result \eqref{distr}.
		\end{theorem}

		\begin{remark}
			When $\nu=1$, equation \eqref{distr} easily reduces to the distribution
			of the classical linear death process, i.e.
			\begin{equation}
				p_k(t) = \binom{n_0}{k} e^{-k \mu t} \left( 1-e^{- \mu t} \right)^{n_0-k},
				\qquad t>0, \: 0 \leq k \leq n_0.
			\end{equation}
		\end{remark}

		In the following theorem we give a proof of an interesting subordination relation.
		\begin{theorem}
			The fractional linear death process $M^\nu (t)$, $t>0$ can be
			represented as
			\begin{equation}
				\label{rec}
				M^\nu(t) \overset{i.d.}{=} M (T_{2\nu}(t)), \qquad t>0, \: \nu \in (0,1],
			\end{equation}
			where $M(t)$, $t>0$ is the classical linear death process
			(see e.g.\ \ocite{bailey}, page 90) and $T_{2 \nu}(t)$, $t>0$, is a random
			time process whose one-dimensional distribution coincides with
			the folded solution to the following fractional diffusion equation
			\begin{equation}
				\label{diffusion}
				\begin{cases}
					\frac{\partial^{2\nu}}{\partial t^{2\nu}} q(s,t) = \frac{\partial^2}{
					\partial s^2} q(s,t), & t>0, \: \nu \in (0,1], \\
					q(s,0) = \delta(s),
				\end{cases}
			\end{equation}
			with the additional condition $q_t(s,0) = 0$ if $\nu \in (1/2,1]$
			(see \ocite{ob}).
			\begin{proof}
				By evaluating the Laplace transform of the generating function of
				the fractional linear death process $M^\nu(t)$, $t>0$, we obtain that
				\begin{align}
					& \int_0^\infty e^{-zt} G^\nu (u,t) dt \\
					& = \int_0^\infty e^{-zt} \sum_{k=0}^{n_0} u^k \binom{n_0}{k}
					\sum_{r=0}^{n_0-k} \binom{n_0-k}{r}
					(-1)^r E_{\nu,1} (-(k+r) \mu t^\nu) dt \notag \\
					& = \sum_0^{n_0} u^k \binom{n_0}{k} \sum_{r=0}^{n_0-k}
					\binom{n_0-k}{r} (-1)^r
					\frac{z^{\nu-1}}{z^\nu+(k+r)\mu} \notag \\
					& = \int_0^\infty \sum_{k=0}^{n_0} u^k \binom{n_0}{k}
					\sum_{r=0}^{n_0-k} \binom{n_0-k}{r} (-1)^r z^{\nu-1}
					e^{-s(z^\nu + (k+r) \mu )} ds \notag \\
					& = \int_0^\infty e^{-sz^\nu} z^{\nu-1} \left\{
					\sum_{k=0}^{n_0} u^k \binom{n_0}{k} \sum_{r=0}^{n_0-k}
					\binom{n_0-k}{r} (-1)^r e^{-s (k+r) \mu} \right\} ds \notag \\
					& = \int_0^\infty e^{-sz^\nu} z^{\nu-1} \left\{
					\sum_{k=0}^{n_0} u^k \binom{n_0}{k} e^{-\mu sk} \sum_{r=0}^{n_0-k}
					\binom{n_0-k}{r} (-1)^r e^{-sr \mu} \right\} ds \notag \\
					& = \int_0^\infty e^{-sz^\nu} z^{\nu-1} \left\{
					\sum_{k=0}^{n_0} u^k \binom{n_0}{k} e^{-\mu sk} (1-e^{-\mu s})^{
					n_0-k} \right\} ds \notag \\
					& = \int_0^\infty e^{-sz^\nu} z^{\nu-1} G(u,s) ds \notag \\
					& = \int_0^\infty e^{-zt} \int_0^\infty \sum_{k=0}^{n_0}
					u^k \text{Pr} \left\{ M(s) = k \right\} f_{T_{2\nu}} (s,t)
					ds \, dt \notag \\
					& = \int_0^\infty e^{-zt} \left\{ \sum_{k=0}^\infty u^k
					\text{Pr} \left\{ M(T_{2 \nu} (t)) = k \right\} \right\} dt, \notag
				\end{align}
				and this is sufficient to prove that \eqref{rec} holds.
				Note that we used two facts. The
				first one is that
				\begin{equation}
					\int_0^\infty e^{-zt} f_{T_{2 \nu}} (s,t) dt = z^{\nu-1} e^{-s z^\nu},
					\qquad s>0, \: z>0,
				\end{equation}
				is the Laplace transform of the solution to \eqref{diffusion}. The second
				fact is that the Laplace transform of the Mittag-Leffler function is
				\begin{equation}
					\label{lap}
					\int_0^\infty e^{-zt} E_{\nu,1} (- \vartheta t^\nu) dt =
					\frac{z^{\nu-1}}{z^\nu + \vartheta}.
				\end{equation}
			\end{proof}
		\end{theorem}
		In figures \ref{nzero} and \ref{nzeromu}, we compare the behaviour of the fractional
		probabilities $p_{n_0}^{0.7} (t)$ and $p_{n_0-1}^{0.7} (t)$ with their classical
		counterparts $p_{n_0}^1 (t)$ and $p_{n_0-1}^1 (t)$, $t>0$. What emerges from the inspection
		of both figures is that, for large values of $t$, the probabilities,
		in the fractional case, decrease more slowly
		than $p_{n_0}^1 (t)$ and $p_{n_0-1}^1 (t)$. The
		probability $p_{n_0-1}^{0.7} (t)$, increases initially faster than $p_{n_0-1}^1 (t)$,
		but after a certain time lapse, $p_{n_0-1}^1 (t)$ dominates $p_{n_0-1}^{0.7} (t)$.

		\begin{remark}
			For $\nu=1/2$, in view of the integral representation
			\begin{equation}
				E_{\frac{1}{2},1} (x) = \frac{2}{\sqrt{\pi}} \int_0^\infty e^{-w^2 +2xw} dw,
				\qquad x \in \mathbb{R},
			\end{equation}
			we extract from \eqref{distr} that
			\begin{align}
				p_k^{\frac{1}{2}}(t) & = \frac{2}{\sqrt{\pi}} \int_0^\infty
				e^{-w^2} \binom{n_0}{k} \sum_{r=0}^{n_0-k} \binom{n_0-k}{r} (-1)^r
				e^{-2w (k+r) \mu t^{\frac{1}{2}}} \\
				& = \int_0^\infty \frac{ e^{-\frac{y^2}{4t}} }{\sqrt{\pi t}} p_k^1 (y) dy =
				\text{Pr} \left\{ M(\left| B(t) \right| ) = k \right\}, \notag
			\end{align}
			where $B(t)$, $t>0$ is a Brownian motion with volatility equal to 2.
		\end{remark}

		\begin{remark}
			We can interpret formula \eqref{distr}
			in an alternative way, as follows. For each integer
			$k \in [0,n_0]$ we have that
			\begin{align}
				p_k^\nu (t) & = \text{Pr} \left\{ M^\nu (t) = k \mid M^\nu (0) = n_0 \right\} \\
				& = \int_0^\infty p_k(s) \text{Pr} \left\{T_{2 \nu} (t) \in ds \right\}
				\notag \\
				& = \binom{n_0}{k} \sum_{r=0}^{n_0-k} \binom{n_0-k}{r} (-1)^r
				\int_0^\infty e^{-\mu (k+r) s} \text{Pr} \left\{ T_{2 \nu} (t) \in ds \right\}
				\notag \\
				& = \binom{n_0}{k} \sum_{r=0}^{n_0-k} \binom{n_0-k}{r} (-1)^r
				\int_0^\infty e^{-\mu (k+r) s} t^{-\nu} W_{-\nu,1-\nu} (-st^{-\nu}) ds \notag \\
				& = \binom{n_0}{k} \sum_{r=0}^{n_0-k} \binom{n_0-k}{r} (-1)^r
				\int_0^\infty e^{- \xi \mu (k+r) t^\nu} W_{-\nu,1-\nu} (-\xi) d \xi \notag \\
				& = \int_0^\infty W_{-\nu,1-\nu} (-\xi) \text{Pr} \left\{
				M_{\xi} (t^\nu) = k \mid M_{\xi} (0) = n_0 \right\} d \xi, \notag
			\end{align}
			where $W_{-\nu, 1-\nu}(-\xi)$ is the Wright function defined as
			\begin{equation}
				W_{-\nu, 1-\nu}(-\xi) = \sum_{r=0}^\infty \frac{(- \xi)^r}{r!
				\Gamma \left( 1- \nu (r+1) \right)}, \qquad 0 < \nu \leq 1.
			\end{equation}
			We therefore obtain an interpretation in terms of a classical linear death
			process $M_{\Xi} (t)$, $t>0$ evaluated on a new time scale and with random
			rate $\mu \cdot \Xi$, where $\Xi$ is a random variable, $\xi \in \mathbb{R}^+$,
			with Wright density
			\begin{equation}
				f_{\Xi} (\xi) = W_{-\nu,1-\nu} (-\xi), \qquad \xi \in \mathbb{R}^+.
			\end{equation}
		\end{remark}

		From equation \eqref{cauchy-fnl} with $\mu_k=k \cdot \mu$,
		the related
		fractional differential equation governing the probability generating function,
		can be easily obtained, leading to
		\begin{equation}
			\begin{cases}
				\frac{\partial^\nu}{\partial t^\nu} G^\nu(u,t) = -\mu u(u-1)
				\frac{\partial}{\partial u} G^\nu(u,t), & \nu \in (0,1], \\
				G^\nu(u,0) = u^{n_0}.
			\end{cases}
		\end{equation}
		From this, and by considering that $\mathbb{E} M^\nu (t) = \left.
		\frac{\partial}{\partial u} G^\nu (u,t)\right|_{u=1}$, we obtain that
		\begin{equation}
			\label{mean}
			\begin{cases}
				\frac{d^\nu}{dt^\nu} \mathbb{E} M^\nu (t) = - \mu \mathbb{E} M^\nu (t),
				& \nu \in (0,1], \\
				\mathbb{E} M^\nu (t) = n_0.
			\end{cases}
		\end{equation}
		Equation \eqref{mean} is easily solved by means of the Laplace transforms, yielding
		\begin{equation}
			\mathbb{E} M^\nu (t) = n_0 E_{\nu,1} (-\mu t^\nu), \qquad t>0, \: \nu \in (0,1].
		\end{equation}

		\begin{remark}
			The mean value $\mathbb{E} M^\nu (t)$ can also be directly
			calculated.
			\begin{align}
				\label{meanlin}
				\mathbb{E} M^\nu (t) & = \sum_{k=0}^{n_0} k p_k^\nu (t) \\
				& = \sum_{k=0}^{n_0} k \binom{n_0}{k} \sum_{r=k}^{n_0}
				\binom{n_0-k}{r-k} (-1)^{r-k} E_{\nu,1} (-r \mu t^\nu) \notag \\
				& = \sum_{r=0}^{n_0} E_{\nu,1} (-r \mu t^\nu ) (-1)^r \sum_{k=1}^r
				k \binom{n_0}{k} \binom{n_0-k}{r-k} (-1)^k \notag \\
				& = \sum_{r=1}^{n_0} E_{\nu,1} (-r \mu t^\nu ) (-1)^r
				n_0 \binom{n_0-1}{r-1} \sum_{k=1}^r \binom{r-1}{k-1} (-1)^k \notag \\
				& = n_0 E_{\nu,1} (- \mu t^\nu). \notag
			\end{align}
			This last step in \eqref{meanlin} holds because
			\begin{align}
				\sum_{k=1}^r \binom{r-1}{k-1} (-1)^k =
				\sum_{k=0}^{r-1} \binom{r-1}{k} (-1)^{k+1} =
				\begin{cases}
					-1, & r=1, \\
					0, & r>1.
				\end{cases}
			\end{align}
		\end{remark}

	\section{Related models}

		\label{gen}
		In this section we present two models which are related to the fractional linear
		death process. The first one is its natural generalisation to the non-linear
		case i.e.\ we consider death rates in the form $\mu_k > 0$, $0 \leq k \leq n_0$.
		The second one is a sublinear process (see \ocite{donnelly}), namely with death rates in the form
		$\mu_k = \mu (n_0+1-k)$; the death rates are thus an increasing sequence as the number
		of individuals in the population decreases towards zero.

		\subsection{Generalisation to the non-linear case}

			Let us denote by $\mathpzc{M}^\nu (t)$, $t>0$ the random number of components
			of a non-linear fractional death process with death rates $\mu_k > 0$,
			$0 \leq k \leq n_0$.

			The state probabilities $\mathpzc{p}_k^\nu (t) = \text{Pr} \left\{ \mathpzc{M}^\nu (t)
			= k \mid \mathpzc{M}^\nu (0) = n_0 \right\}$, $t>0$, $0 \leq k \leq n_0$,
			$\nu \in (0,1]$ are governed by the following difference-differential equations
			\begin{equation}
				\label{cn}
				\begin{cases}
					\frac{d^\nu}{dt^\nu} \mathpzc{p}_k(t) = \mu_{k+1}
					\mathpzc{p}_{k+1} (t) - \mu_k \mathpzc{p}_k (t),
					& 0 < k < n_0, \\
					\frac{d^\nu}{dt^\nu} \mathpzc{p}_0 (t) = \mu_1 \mathpzc{p}_1 (t),
					& k=0, \\
					\frac{d^\nu}{dt^\nu} \mathpzc{p}_{n_0}(t) = -\mu_{n_0}
					\mathpzc{p}_{n_0} (t), & k = n_0, \\
					\mathpzc{p}_k (0) =
						\begin{cases}
							0, & 0 \leq k < n_0, \\
							1, & k=n_0.
						\end{cases}
				\end{cases}
			\end{equation}

			The fractional derivatives appearing in \eqref{cn} provide the system with a
			global memory; i.e.\ the evolution of the state probabilities $\mathpzc{p}_k^\nu (t)$,
			$t>0$, is influenced by the past, as definition \eqref{caputo} shows.
			This is a major difference with the classical non-linear (and, of course,
			linear and sublinear) death processes, and reverberates in the slowly
			decaying structure of probabilities extracted from \eqref{cn}.

			In the non-linear process, the dependence of death rates from the size of the population
			is arbitrary, and this explains the complicated structure of the
			probabilities obtained.
			Further generalisation can be considered by assuming that the death rates
			depend on $t$ (non-homogeneous, non-linear death process).

			We outline here the evaluation of the probabilities $\mathpzc{p}_k^\nu (t)$,
			$t > 0$, $0 \leq k \leq n_0$, which can be obtained, as in the linear case,
			by means of a recursive procedure (similar to that implemented in
			\ocite{pol} for the fractional linear birth process).

			Let $k = n_0$. By means of the Laplace transform applied to equation \eqref{cn} we
			immediately obtain that
			\begin{equation}
				\mathpzc{p}_{n_0}^\nu (t) = E_{\nu,1} (- \mu_{n_0} t^\nu).
			\end{equation}
			When $k=n_0-1$ we get
			\begin{align}
				& z^\nu \mathcal{L} \left\{ \mathpzc{p}_{n_0-1}^\nu \right\} (z) = - \mu_{n_0-1}
				\mathcal{L} \left\{ \mathpzc{p}_{n_0-1}^\nu \right\} (z)
				+ \mu_{n_0} \frac{z^{\nu-1}}{
				z^\nu + \mu_{n_0}} \\
				& \Leftrightarrow \mathcal{L} \left\{ \mathpzc{p}_{n_0-1}^\nu \right\} (z) =
				\mu_{n_0} \frac{z^{\nu-1}}{z^\nu + \mu_{n_0}} \cdot \frac{1}{z^\nu + \mu_{n_0-1}}
				\notag \\
				& \Leftrightarrow \mathcal{L} \left\{ \mathpzc{p}_{n_0-1}^\nu \right\} (z) =
				\mu_{n_0} z^{\nu-1} \left[ \frac{1}{z^\nu + \mu_{n_0}}
				- \frac{1}{z^\nu+\mu_{n_0-1}}
				\right] \frac{1}{\mu_{n_0-1} - \mu_{n_0}} \notag \\
				& \Leftrightarrow \mathpzc{p}_{n_0-1}^\nu (t) =
				\frac{\mu_{n_0}}{\mu_{n_0-1} - \mu_{n_0}} \bigg\{ E_{\nu,1} (- \mu_{n_0} t^\nu )
				- E_{\nu,1} (- \mu_{n_0-1} t^\nu ) \bigg\}. \notag
			\end{align}

			For $k=n_0-2$ we obtain in the same way that
			\begin{align}
				& z^\nu \mathcal{L} \left\{ \mathpzc{p}_{n_0-2}^\nu \right\} (z) \\
				& = - \mu_{n_0-2} \mathcal{L} \left\{ \mathpzc{p}_{n_0-2}^\nu \right\} (z)
				+ \frac{\mu_{n_0} \mu_{n_0-1}}{\mu_{n_0-1}-\mu_{n_0}} \left[
				\frac{z^{\nu-1}}{z^\nu+\mu_{n_0}} - \frac{z^{\nu-1}}{z^\nu +
				\mu_{n_0-1}} \right],
				\notag
			\end{align}
			so that
			\begin{align}
				& \hspace{-0.5cm} \mathcal{L} \left\{  \mathpzc{p}_{n_0-2}^\nu \right\} (z) \\
				= {} &
				\frac{\mu_{n_0} \mu_{n_0-1}}{\mu_{n_0-1}-\mu_{n_0}} z^{\nu-1}
				\left[ \frac{1}{z^\nu + \mu_{n_0}} - \frac{1}{z^\nu+\mu_{n_0-1}} \right]
				\frac{1}{z^\nu+\mu_{n_0-2}} \notag \\
				= {} & \frac{\mu_{n_0} \mu_{n_0-1}}{\mu_{n_0-1}-\mu_{n_0}} z^{\nu-1}
				\left[ \left( \frac{1}{z^\nu+\mu_{n_0}} - \frac{1}{z^\nu + \mu_{n_0-2}} \right)
				\frac{1}{\mu_{n_0-2} - \mu_{n_0}} \right. \notag \\
				& \left. - \left( \frac{1}{z^\nu + \mu_{n_0-1}}
				- \frac{1}{z^\nu+\mu_{n_0-2}} \right)
				\frac{1}{\mu_{n_0-2} - \mu_{n_0-1}} \right]. \notag
			\end{align}
			By inverting the Laplace transform we readily arrive at the following result
			\begin{align}
				\mathpzc{p}_{n_0-2}^\nu {} & (t) = \mu_{n_0} \mu_{n_0-1} \left[
				\frac{E_{\nu,1} (-
				\mu_{n_0} t^\nu )}{(\mu_{n_0-1}-\mu_{n_0})(\mu_{n_0-2} -\mu_{n_0})} \right. \\
				& - \frac{E_{\nu,1} (-\mu_{n_0-2}
				t^\nu )}{(\mu_{n_0-1}-\mu_{n_0})(\mu_{n_0-2} -\mu_{n_0})}
				- \frac{E_{\nu,1} (-\mu_{n_0-1}
				t^\nu )}{(\mu_{n_0-1}-\mu_{n_0})(\mu_{n_0-2} -\mu_{n_0-1})} \notag \\
				& \left. + \frac{E_{\nu,1} (-\mu_{n_0-2}
				t^\nu )}{(\mu_{n_0-1}-\mu_{n_0})(\mu_{n_0-2} -\mu_{n_0-1})} \right] \notag \\
				= {} & \mu_{n_0} \mu_{n_0-1} \left[ \frac{E_{\nu,1} (-
				\mu_{n_0} t^\nu )}{(\mu_{n_0-1}-\mu_{n_0})(\mu_{n_0-2}
				-\mu_{n_0})} \right. \notag \\
				& + \frac{E_{\nu,1} (-\mu_{n_0-2} t^\nu)}{(\mu_{n_0-1}-\mu_{n_0})} \left(
				\frac{1}{\mu_{n_0-2}-\mu_{n_0-1}}
				- \frac{1}{\mu_{n_0-2}-\mu_{n_0}} \right) \notag \\
				& \left. - \frac{E_{\nu,1} (- \mu_{n_0-1}
				t^\nu )}{(\mu_{n_0-1}-\mu_{n_0})(\mu_{n_0-2}
				-\mu_{n_0-1})} \right] \notag \\
				= {} & \mu_{n_0} \mu_{n_0-1} \left[ \frac{E_{\nu,1} (-
				\mu_{n_0} t^\nu )}{(\mu_{n_0-1}-\mu_{n_0})(\mu_{n_0-2}
				-\mu_{n_0})} \right. \notag \\
				& \left. + \frac{E_{\nu,1} (-\mu_{n_0-2} t^\nu )}{(\mu_{n_0-2}-
				\mu_{n_0-1})(\mu_{n_0-2} -\mu_{n_0})}
				- \frac{E_{\nu,1} (-\mu_{n_0-1} t^\nu
				)}{(\mu_{n_0-1}-\mu_{n_0})(\mu_{n_0-2} -\mu_{n_0-1})}
				\right]. \notag
			\end{align}
			The structure of the state
			probabilities for arbitrary values of $k= n_0-l$, $0 \leq l < n_0$,
			can now be easily obtained.
			The proof follows the lines of the derivation of the
			state probabilities for the fractional non-linear pure birth process
			adopted in Theorem 2.1 in \ocite{pol}. We have that
			\begin{equation}
				\label{nlinearnu}
				\mathpzc{p}_{n_0-l}^\nu (t) =
				\begin{cases}
					\textstyle\prod\limits_{j=0}^{l-1} \mu_{n_0-j}
					\displaystyle\sum_{m=0}^l \frac{ E_{\nu,1} (-
					\mu_{n_0-m} t^\nu)}{\textstyle\prod\limits_{\substack{
					h=0 \\ h \neq m}}^l \left( \mu_{n_0-h} -
					\mu_{n_0-m} \right) }, & 1 \leq l < n_0, \\
					E_{\nu,1} (-\mu_{n_0} t^\nu), & l=0.
				\end{cases}
			\end{equation}
			By means of some changes of indices, formula \eqref{nlinearnu} can also
			be written as
			\begin{equation}
				\label{nlinearnu2}
				\mathpzc{p}_k^\nu (t) =
				\begin{cases}
					\textstyle\prod\limits_{j=k+1}^{n_0} \mu_j
					\displaystyle\sum_{m=k}^{n_0} \frac{ E_{\nu,1} (-
					\mu_m t^\nu)}{\textstyle\prod\limits_{\substack{
					h=k \\ h \neq m}}^{n_0} \left( \mu_h -
					\mu_m \right) }, & 0 < k < n_0, \\
					E_{\nu,1} (-\mu_{n_0} t^\nu), & k=n_0.
				\end{cases}
			\end{equation}

			For the extinction probability, we have to solve the following initial value problem:
			\begin{equation}
				\begin{cases}
					\label{zz}
					\frac{d^\nu}{dt^\nu} \mathpzc{p}_0 (t) = \mu_1
					\textstyle\prod\limits_{j=2}^{n_0} \mu_j
					\displaystyle\sum_{m=1}^{n_0} \frac{ E_{\nu,1} (-
					\mu_m t^\nu)}{\textstyle\prod\limits_{\substack{
					h=1 \\ h \neq m}}^{n_0} \left( \mu_h -
					\mu_m \right) }, & n_0 > 1, \\
					\frac{d^\nu}{dt^\nu} \mathpzc{p}_0 (t) = \mu_1
					E_{\nu,1} (- \mu_1 t^\nu), & n_0=1, \\
					\mathpzc{p}_0 (0) = 0, & n_0 \geq 1.
				\end{cases}
			\end{equation}
			When $n_0>1$, starting from \eqref{zz} and by resorting to the Laplace transform
			once again, we have that
			\begin{equation}
				\label{inv}
				\mathcal{L} \left\{ \mathpzc{p}_0^\nu \right\} (z) =
				\textstyle\prod\limits_{j=1}^{n_0} \mu_j
				\displaystyle\sum_{m=1}^{n_0}
				\frac{1}{\textstyle\prod\limits_{\substack{
				h=1 \\ h \neq m}}^{n_0} \left( \mu_h -
				\mu_m \right) } \cdot \frac{z^{-1}}{z^\nu+\mu_m}.
			\end{equation}
			The inverse Laplace transform of \eqref{inv} leads to
			\begin{align}
				\mathpzc{p}_0^\nu (t) & =
				\textstyle\prod\limits_{j=1}^{n_0} \mu_j
				\displaystyle\sum_{m=1}^{n_0}
				\frac{1}{\textstyle\prod\limits_{\substack{
				h=1 \\ h \neq m}}^{n_0} \left( \mu_h -
				\mu_m \right) } t^\nu E_{\nu,\nu+1} (-\mu_m t^\nu ) \\
				& = \textstyle\prod\limits_{j=1}^{n_0} \mu_j
				\displaystyle\sum_{m=1}^{n_0}
				\frac{1}{\textstyle\prod\limits_{\substack{
				h=1 \\ h \neq m}}^{n_0} \left( \mu_h -
				\mu_m \right) } \cdot \frac{1}{\mu_m} \left[ 1- E_{\nu,1} (-\mu_m t^\nu)
				\right] \notag \\
				& = \sum_{m=1}^{n_0} \prod_{\substack{
				h=1 \\ h \neq m}}^{n_0} \left( \frac{\mu_h}{\mu_h-\mu_m} \right)
				- \sum_{m=1}^{n_0} \prod_{\substack{
				h=1 \\ h \neq m}}^{n_0} \left( \frac{\mu_h}{\mu_h-\mu_m} \right)
				E_{\nu,1} (- \mu_m t^\nu ) \notag \\
				& = 1 - \sum_{m=1}^{n_0} \prod_{\substack{
				h=1 \\ h \neq m}}^{n_0} \left( \frac{\mu_h}{\mu_h-\mu_m} \right)
				E_{\nu,1} (- \mu_m t^\nu ). \notag
			\end{align}
			Note that, in the last step, we used the following fact:
			\begin{equation}
				\sum_{m=1}^{n_0} \prod_{\substack{
				h=1 \\ h \neq m}}^{n_0} \left( \frac{\mu_h}{\mu_h-\mu_m} \right)
				= 1.
			\end{equation}
			This can be ascertained by observing that
			\begin{equation}
				\prod_{1 \leq h < l \leq n_0} (\mu_h - \mu_l) =
				\det \bm{A}
				= \sum_{j=1}^{n_0} a_{1,j} (-1)^{j+1} \text{Min}_{1,j}
			\end{equation}
			where
			\begin{equation}
				\bm{A} = \left|
				\begin{array}{cccc}
					1 & 1 & \dots & 1 \\
					\mu_1 & \mu_2 & \dots & \mu_{n_0} \\
					\mu_1^2 & \mu_2^2 & \dots & \mu_{n_0}^2 \\
					\vdots & \vdots & \ddots & \vdots \\
					\mu_1^{n_0-1} & \mu_2^{n_0-1} & \dots & \mu_{n_0}^{n_0-1}
				\end{array}
				\right|,
			\end{equation}
			is a Vandermonde matrix and $\text{Min}_{1,j}$ is the determinant of the
			matrix resulting from $\bm{A}$ by removing the first row and the
			$j$-th column.

			When $n_0=1$ we obtain
			\begin{equation}
				\mathcal{L} \left\{ \mathpzc{p}_0^\nu \right\} (z) =
				\mu_1 \frac{z^{-1}}{z^\nu + \mu_1},
			\end{equation}
			so that the inverse Laplace transform can be written as
			\begin{align}
				\mathpzc{p}_0^\nu (t) & = \mu_1 t^\nu E_{\nu,\nu+1} (- \mu_1 t^\nu ) \\
				& = 1 - E_{\nu,1} ( - \mu_1 t^\nu ). \notag
			\end{align}

			We can therefore summarise the results obtained as follows:
			\begin{equation}
				\label{nlinearnu3}
				\mathpzc{p}_k^\nu (t) =
				\begin{cases}{}
					\textstyle\prod\limits_{j=k+1}^{n_0} \mu_j
					\displaystyle\sum_{m=k}^{n_0} \frac{ E_{\nu,1} (-
					\mu_m t^\nu)}{\textstyle\prod\limits_{\substack{
					h=k \\ h \neq m}}^{n_0} \left( \mu_h -
					\mu_m \right) }, & 0 < k < n_0, \: n_0 > 1, \\
					E_{\nu,1} (-\mu_{n_0} t^\nu), & k=n_0, \: n_0 \geq 1,
				\end{cases}
			\end{equation}
			and
			\begin{equation}
				\mathpzc{p}_0^\nu (t) =
				\begin{cases}
					1 - \displaystyle\sum_{m=1}^{n_0} \displaystyle\prod_{\substack{
					h=1 \\ h \neq m}}^{n_0} \left( \frac{\mu_h}{\mu_h-\mu_m} \right)
					E_{\nu,1} (- \mu_m t^\nu ), & n_0 > 1, \\
					1 - E_{\nu,1} (- \mu_1 t^\nu), & n_0 = 1.
				\end{cases}
			\end{equation}

		\subsection{A fractional sublinear death process}

			\label{alternative}
			We consider in this section the process where the infinitesimal death probabilities
			have the form
			\begin{equation}
				\text{Pr} \left\{ \mathfrak{M} (t,t+dt] = -1 \mid
				\mathfrak{M} (t) = k \right\} = \mu (n_0 +1 -k) dt + o(dt),
			\end{equation}
			where $n_0$ is the initial number of individuals in the population.
			The state probabilities
			\begin{equation}
				\mathfrak{p}_k (t) = \text{Pr} \left\{ \mathfrak{M} (t) = k \mid
				\mathfrak{M} (0) = n_0 \right\}, \qquad 0 \leq k \leq n_0,
			\end{equation}
			satisfy the equations
			\begin{equation}
				\begin{cases}
					\frac{d}{dt} \mathfrak{p}_k(t) =
					- \mu (n_0+1-k) \mathfrak{p}_k (t)
					+ \mu (n_0-k) \mathfrak{p}_{k+1}(t),
					& 1 \leq k \leq n_0, \\
					\frac{d}{dt} \mathfrak{p}_0 (t) = \mu n_0 \mathfrak{p}_1 (t),
					& k = 0,  \\
					\mathfrak{p}_k (0) =
						\begin{cases}
							1, & k=n_0, \\
							0, & 0 \leq k < n_0.
						\end{cases}
				\end{cases}
			\end{equation}
			In this model the death rate increases with decreasing population
			size.

			The probabilities
			$\mathfrak{p}_k^\nu (t) = \text{Pr} \left\{ \mathfrak{M}^\nu (t) = k
			\mid \mathfrak{M}^\nu (0) = n_0 \right\}$ of the
			fractional version of this process are governed by the equations
			\begin{equation}
				\label{cna}
				\begin{cases}
					\frac{d^\nu}{dt^\nu} \mathfrak{p}_k(t) =
					- \mu (n_0+1-k) \mathfrak{p}_k (t)
					+ \mu (n_0-k) \mathfrak{p}_{k+1}(t),
					& 1 \leq k \leq n_0, \\
					\frac{d^\nu}{dt^\nu} \mathfrak{p}_0 (t) = \mu n_0 \mathfrak{p}_1 (t),
					& k = 0,  \\
					\mathfrak{p}_k (0) =
						\begin{cases}
							1, & k=n_0, \\
							0, & 0 \leq k < n_0.
						\end{cases}
				\end{cases}
			\end{equation}

			We first observe that the solution to the Cauchy problem
			\begin{equation}
				\begin{cases}
					\frac{d^\nu}{dt^\nu} \mathfrak{p}_{n_0} (t)
					= - \mu \mathfrak{p}_{n_0} (t), \\
					\mathfrak{p}_{n_0} (0) = 1,
				\end{cases}
			\end{equation}
			is $\mathfrak{p}_{n_0}^\nu (t) = E_{\nu,1} (- \mu t^\nu)$, $t>0$.
			\begin{figure}
				\centering
				\includegraphics[scale=0.8,angle=0]{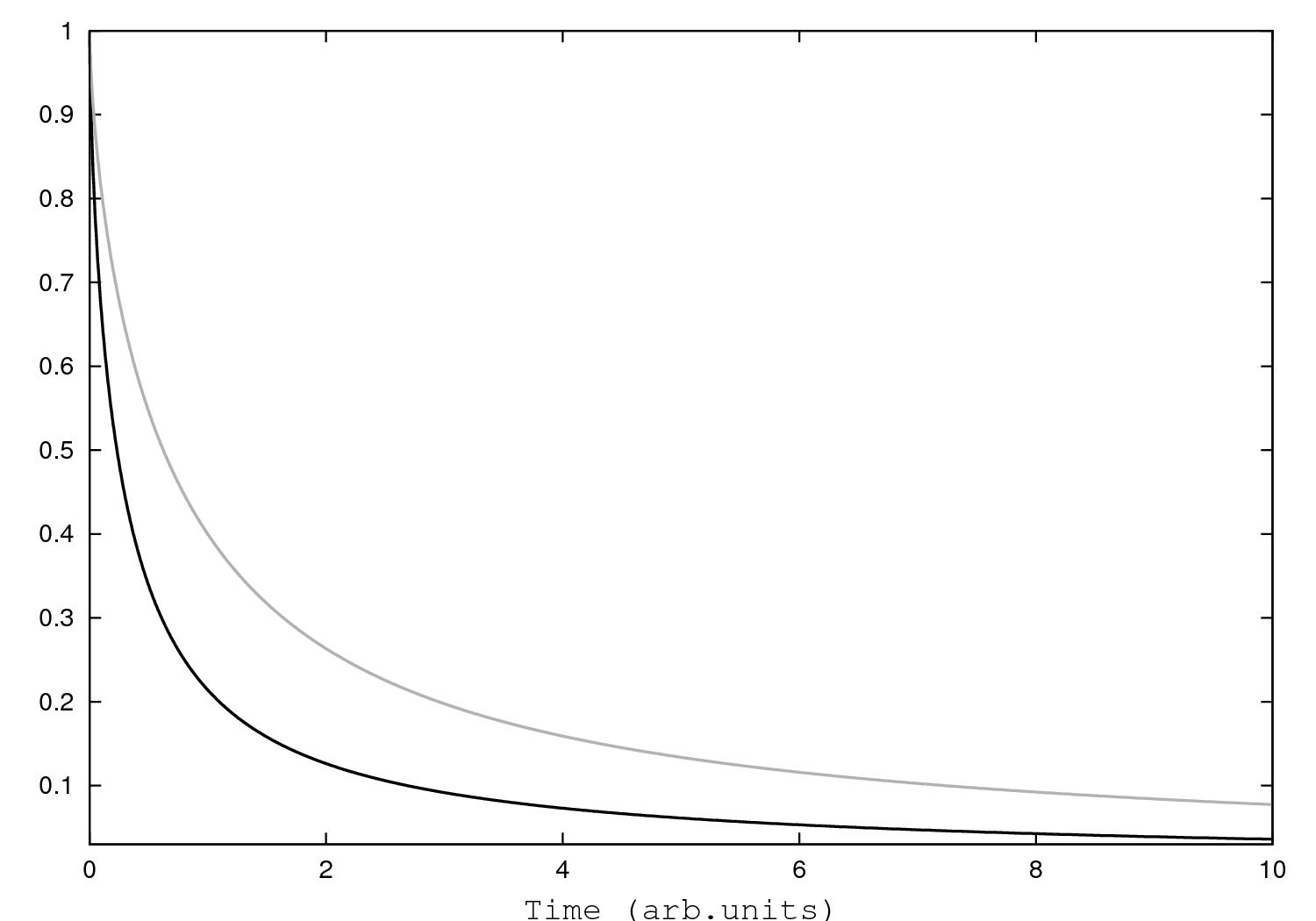}
				\caption{\label{pplot2}Plot of $p_{n_0}^{0.7}
					(t)$ (in black) and $\mathfrak{p}_{n_0}^{0.7}
					(t)$ (in grey), with $n_0 = 2$.}
			\end{figure}

			In order to solve the equation
			\begin{equation}
				\begin{cases}
					\frac{d^\nu}{dt^\nu} \mathfrak{p}_{n_0-1} (t)
					= - 2 \mu \mathfrak{p}_{n_0-1} (t)
					+ \mu E_{\nu,1} (- \mu t^\nu), \\
					\mathfrak{p}_{n_0-1} (0) = 0,
				\end{cases}
			\end{equation}
			we resort to the Laplace transform and obtain that
			\begin{align}
				\label{z}
				\mathcal{L} \left\{ \mathfrak{p}_{n_0-1}^\nu \right\}(z) & =
				\mu z^{\nu -1} \frac{1}{z^\nu +\mu} \cdot \frac{1}{
				z^\nu + 2 \mu} \\
				& = z^{\nu-1} \left( \frac{1}{z^\nu + \mu} -
				\frac{1}{z^\nu + 2 \mu} \right). \notag
			\end{align}
			By inverting \eqref{z} we extract the following result
			\begin{equation}
				\mathfrak{p}_{n_0-1}^\nu (t) = E_{\nu,1} (-\mu t^\nu) -
				E_{\nu,1} (-2 \mu t^\nu ).
			\end{equation}

			By the same technique we solve
			\begin{equation}
				\begin{cases}
					\frac{d^\nu}{dt^\nu} \mathfrak{p}_{n_0-2} (t)
					= - 3 \mu \mathfrak{p}_{n_0-2} (t)
					+ 2\mu
					\left[ E_{\nu,1} (- \mu t^\nu)
					- E_{\nu,1} (-2 \mu t^\nu) \right], \\
					\mathfrak{p}_{n_0-2} (0) = 0,
				\end{cases}
			\end{equation}
			thus obtaining
			\begin{align}
				\label{w}
				\mathcal{L} & \left\{ \mathfrak{p}_{n_0-2}^\nu \right\} (z) =
				2 \mu z^{\nu-1}
				\left[ \frac{1}{z^\nu + \mu} - \frac{1}{z^\nu+2\mu} \right]
				\frac{1}{z^\nu+3\mu} \\
				& = 2\mu z^{\nu-1}
				\left[ \left( \frac{1}{z^\nu+\mu} - \frac{1}{z^\nu + 3\mu} \right)
				\frac{1}{2\mu}
				- \left( \frac{1}{z^\nu + 2\mu}
				- \frac{1}{z^\nu+3\mu} \right)
				\frac{1}{\mu} \right] \notag \\
				& = \frac{z^{\nu-1}}{z^\nu+\mu} - 2 \frac{z^{\nu-1}}{z^\nu+2 \mu}
				+ \frac{z^{\nu-1}}{z^\nu + 3\mu}. \notag
			\end{align}
			In light of \eqref{w}, we infer that
			\begin{equation}
				\mathfrak{p}_{n_0-2}^\nu (t) = E_{\nu,1} (-\mu t^\nu) -2 E_{\nu,1} (-2 \mu t^\nu)
				+ E_{\nu,1} (-3 \mu t^\nu).
			\end{equation}
			\begin{figure}
				\centering
				\includegraphics[scale=0.8,angle=0]{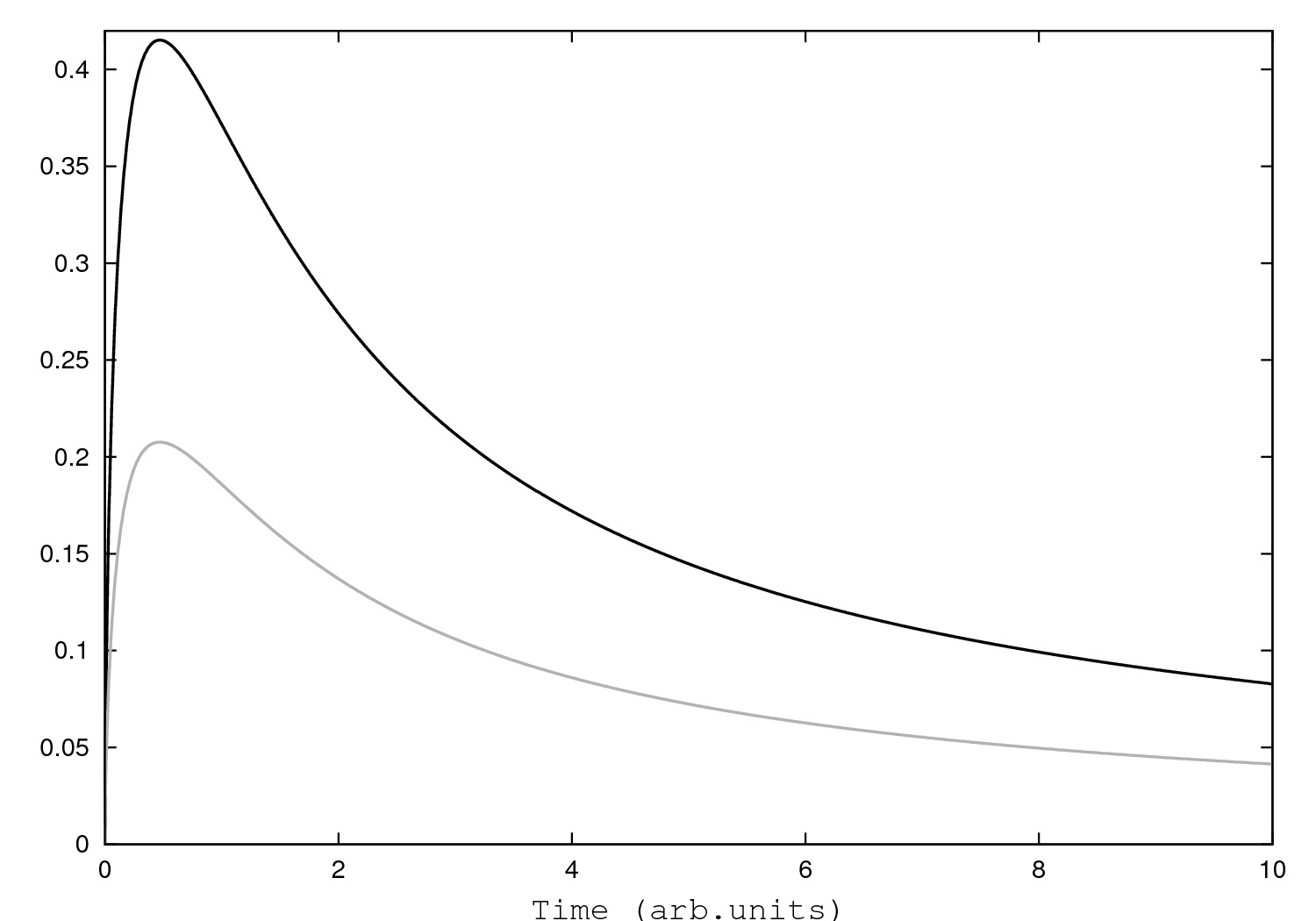}
				\caption{\label{pplot}Plot of $p_{n_0-1}^{0.7}
					(t)$ (in black) and $\mathfrak{p}_{n_0-1}^{0.7}
					(t)$ (in grey), with $n_0 = 2$.}
			\end{figure}

			For all $1 \leq n_0-m \leq n_0$, by similar calculations, we arrive at the general result
			\begin{equation}
				\label{x}
				\mathfrak{p}_{n_0-m}^\nu (t) = \sum_{l=0}^m \binom{m}{l} (-1)^l E_{\nu,1} \left(
				- \left( l+1 \right) \mu t^\nu \right), \qquad 1 \leq n_0-m \leq n_0.
			\end{equation}
			Introducing the notation $n_0-m=k$, we rewrite the state probabilities
			\eqref{x} in the following manner
			\begin{equation}
				\label{pik}
				\mathfrak{p}_k^\nu (t) = \sum_{l=0}^{n_0-k} \binom{n_0-k}{l} (-1)^l
				E_{\nu,1} \left( - \left( l+1 \right)
				\mu t^\nu \right), \qquad 1 \leq k \leq n_0.
			\end{equation}

			For the extinction probability we must solve the following Cauchy problem
			\begin{equation}
				\label{y}
				\begin{cases}
					\frac{d^\nu}{dt^\nu} \mathfrak{p}_0 (t) =
					\mu n_0 \sum_{l=0}^{n_0-1} \binom{n_0-1}{l} (-1)^l
					E_{\nu,1} \left( - \left( l+1 \right)
					\mu t^\nu \right), \\
					\mathfrak{p}_0 (0) = 0.
				\end{cases}
			\end{equation}
			The Laplace transform of \eqref{y} yields
			\begin{equation}
				z^\nu \mathcal{L} \left\{ \mathfrak{p}_0^\nu \right\} (z)
				= \mu n_0 \sum_{l=0}^{n_0-1} \binom{n_0-1}{l} (-1)^l
				\frac{z^{\nu-1}}{z^\nu + \mu (l+1)}.
			\end{equation}
			The inverse Laplace transform can be written down as
			\begin{equation}
				\label{v}
				\mathfrak{p}_0^\nu (t) = \mu n_0 \sum_{l=0}^{n_0-1} \binom{n_0-1}{l}
				(-1)^l \frac{1}{\Gamma (\nu)} \int_0^t E_{\nu,1} \left(-
				(l+1) \mu s^\nu \right) (t-s)^{\nu-1} ds.
			\end{equation}
			The integral appearing in \eqref{v} can be suitably
			evaluated as follows
			\begin{align}
				\label{vv}
				& \int_0^t E_{\nu,1} \left(-
				(l+1) \mu s^\nu \right) (t-s)^{\nu-1} \\
				& = \sum_{m=0}^\infty \frac{(-(l+1)\mu)^m}{\Gamma (\nu m+1)}
				\int_0^t s^{\nu m} (t-s)^{\nu-1} ds \notag \\
				& = \sum_{m=0}^\infty \frac{(-(l+1)\mu)^m}{\Gamma (\nu m+1)}
				\frac{t^{\nu \left( m+1 \right)} \Gamma (\nu) \Gamma ( \nu m+1)}{
				\Gamma (\nu m+\nu +1)} \notag \\
				& = \frac{\Gamma (\nu)}{(-\mu (l+1))} \sum_{m=0}^\infty
				\frac{(-(l+1)\mu t^\nu)^{m+1}}{\Gamma (\nu (m+1)+1)} \notag \\
				& = \frac{\Gamma (\nu)}{(-\mu (l+1))} \left[
				E_{\nu,1} (- (l+1) \mu t^\nu)-1 \right]. \notag
			\end{align}
			By inserting result \eqref{vv} into \eqref{v}, we obtain
			\begin{align}
				\label{extinction}
				\mathfrak{p}_0^\nu (t) & = n_0 \sum_{l=0}^{n_0-1} \binom{n_0-1}{l}
				\frac{(-1)^{l+1}}{
				l+1} \left[ E_{\nu,1} (- (l+1) \mu t^\nu)-1 \right] \\
				& = \sum_{l=0}^{n_0-1} \binom{n_0}{l+1} (-1)^{l+1}
				\left[ E_{\nu,1} (- (l+1) \mu t^\nu)-1 \right] \notag \\
				& = \sum_{l=1}^{n_0} \binom{n_0}{l} (-1)^l E_{\nu,1}
				(-l \mu t^\nu) - \sum_{l=1}^{n_0} \binom{n_0}{l}
				(-1)^l \notag \\
				& = 1 + \sum_{l=1}^{n_0} \binom{n_0}{l} (-1)^l E_{\nu,1} (-l \mu t^\nu)
				\notag \\
				& = \sum_{l=0}^{n_0} \binom{n_0}{l} (-1)^l E_{\nu,1} (-l \mu t^\nu). \notag
			\end{align}

			\begin{remark}
				We check that the probabilities \eqref{pik} and \eqref{extinction}
				sum up to unity. We start by analysing the following sum:
				\begin{equation}
					\label{summa}
					\sum_{k=1}^{n_0} \mathfrak{p}_k^\nu (t) =
					\sum_{k=1}^{n_0} \sum_{l=0}^{n_0-k} \binom{n_0-k}{l}
					(-1)^l E_{\nu,1} (- (l+1) \mu t^\nu ).
				\end{equation}
				In order to evaluate \eqref{summa}, we resort to the Laplace transform
				\begin{equation}
					\sum_{k=1}^{n_0} \mathcal{L} \left\{ \mathfrak{p}_k^\nu \right\}
					(z) = \frac{z^{\nu-1}}{\mu} \sum_{k=1}^{n_0}
					\sum_{l=0}^{n_0-k} \binom{n_0-k}{l} (-1)^l \frac{1}{\frac{z^\nu}{\mu}
					+1+l}.
				\end{equation}
				By using formula (6) of \ocite{kirsch} (see also \ocite{knuth}, formula (5.41),
				page 188),
				we obtain that
				\begin{align}
					\label{check}
					\sum_{k=1}^{n_0} \mathcal{L} \left\{ \mathfrak{p}_k^\nu \right\}
					(z) & = \frac{z^{\nu-1}}{\mu} \sum_{k=1}^{n_0}
					\frac{\Gamma (n_0 -k+1)}{\left( \frac{z^\nu}{\mu}+1 \right)
					\left( \frac{z^\nu}{\mu}+2 \right) \dots \left( \frac{z^\nu}{\mu}
					+1+n_0-k \right)} \\
					& = \frac{z^{\nu-1}}{\mu} \sum_{k=1}^{n_0} \frac{\Gamma \left(
					\frac{z^\nu}{\mu} +1 \right) \Gamma \left( n_0-k+1 \right)}{
					\Gamma \left( \frac{z^\nu}{\mu} +1 + n_0 -k \right)} \notag \\
					& = \frac{z^{\nu-1}}{\mu} \sum_{k=1}^{n_0}
					\int_0^1 x^{\frac{z^\nu}{\mu}} (1-x)^{n_0-k} dx \notag \\
					& = \frac{z^{\nu-1}}{\mu} \int_0^1 x^{\frac{z^\nu}{\mu} -1}
					\left[ 1-(1-x)^{n_0} \right] dx \notag \\
					& = \frac{1}{z} - \frac{z^{\nu-1}}{\mu} \int_0^1
					x^{\frac{z^\nu}{\mu} -1} (1-x)^{n_0} dx \notag \\
					& \overset{\left(- \ln x = y \right)}{=}
					\frac{1}{z} - \frac{z^{\nu-1}}{\mu} \int_0^\infty
					e^{-y \frac{z^\nu}{\mu}} \left( 1- e^{-y} \right)^{n_0} dy \notag \\
					& \overset{\left( y/\mu=w \right)}{=}
					\frac{1}{z} z^{\nu-1} \int_0^\infty e^{-wz^\nu}
					\left( 1-e^{-\mu w} \right)^{n_0} dw \notag \\
					& = \frac{1}{z} - z^{\nu-1} \sum_{k=0}^{n_0}
					\binom{n_0}{k} (-1)^k \int_0^\infty e^{-z^\nu w - \mu w k} dw \notag \\
					& = \frac{1}{z} - z^{\nu-1} \sum_{k=0}^{n_0} \binom{n_0}{k}
					(-1)^k \frac{1}{z^\nu + \mu k}. \notag
				\end{align}
				The inverse Laplace transform of \eqref{check} is therefore
				\begin{align}
					\label{summa2}
					\sum_{k=1}^{n_0} \mathfrak{p}_k^\nu (t) & =
					1 - \sum_{k=0}^{n_0} \binom{n_0}{k} (-1)^k E_{\nu,1} (
					-\mu k t^\nu) \\
					& = - \sum_{k=1}^{n_0} \binom{n_0}{k} (-1)^k E_{\nu,1}
					(-\mu k t^\nu). \notag
				\end{align}
				By putting \eqref{extinction} and \eqref{summa2} together,
				we conclude that
				\begin{equation}
					\sum_{k=0}^{n_0} \mathfrak{p}_k^\nu (t) = 1,
				\end{equation}
				as it should be.
			\end{remark}

			\begin{remark}
				We observe that, in the linear and sublinear
				death processes, the extinction probabilities
				coincide. This implies that although the state probabilities
				$p_k^\nu (t)$ and $\mathfrak{p}_k^\nu (t)$ differ (see figures 3 and 4)
				for all $1 \leq k \leq n_0$,
				we have that
				\begin{equation}
					\sum_{k=1}^{n_0} p_k^\nu (t) = \sum_{k=1}^{n_0} \mathfrak{p}_k^\nu
					(t).
				\end{equation}
				This can be checked by performing the following sum
				\begin{align}
					& \sum_{k=1}^{n_0} \mathcal{L} \left\{ p_k^\nu (t) \right\} (z) \\
					& = \sum_{k=1}^{n_0} \binom{n_0}{k} \sum_{r=0}^{n_0-k}
					\binom{n_0-k}{r} (-1)^r \frac{z^{\nu-1}}{z^\nu + \mu (k+r)} \notag \\
					& = \frac{z^{\nu-1}}{\mu} \sum_{k=1}^{n_0} \binom{n_0}{k}
					\sum_{r=0}^{n_0-k} \binom{n_0-k}{r} (-1)^r \frac{1}{\frac{z^\nu}{\mu}
					+k+r} \notag \\
					& = \frac{z^{\nu-1}}{\mu} \sum_{k=1}^{n_0} \binom{n_0}{k}
					\frac{\left(n_0-k \right)!}{\left( \frac{z^\nu}{\mu} +k \right)
					\left( \frac{z^\nu}{\mu} +k+1 \right) \dots \left(
					\frac{z^\nu}{\mu} +n_0 \right)} \notag \\
					& = \frac{z^{\nu-1}}{\mu} \sum_{k=1}^{n_0} \binom{n_0}{k}
					\frac{ \Gamma \left( n_0-k+1 \right) \Gamma \left( \frac{z^\nu}{\mu}
					+k \right)}{\Gamma \left( \frac{z^\nu}{\mu} +n_0 +1 \right)} \notag \\
					& = \frac{z^{\nu-1}}{\mu} \int_0^1 (1-x)^{\frac{z^\nu}{\mu}-1}
					\sum_{k=1}^{n_0} \binom{n_0}{k} x^{n_0-k} (1-x)^k dx \notag \\
					& = \frac{z^{\nu-1}}{\mu} \int_0^1 (1-x)^{\frac{z^\nu}{\mu}-1}
					\left(1-x^{n_0} \right) dx \notag \\
					& =\frac{1}{z} - \frac{z^{\nu-1}}{\mu} \int_0^1
					x^{n_0} (1-x)^{\frac{z^\nu}{\mu}-1} dx. \notag
				\end{align}
				This coincides with the fourth-to-last step of \eqref{check} and therefore we can
				conclude that
				\begin{equation}
					\sum_{k=1}^{n_0} p_k^\nu (t) = - \sum_{k=1}^{n_0}
					\binom{n_0}{k} (-1)^k E_{\nu,1} (- \mu k t^\nu )
					= \sum_{k=1}^{n_0} \mathfrak{p}_k^\nu (t).
				\end{equation}
			\end{remark}

			\subsubsection{Mean value}

				\begin{theorem}
					Consider the fractional sublinear death process $\mathfrak{M}^\nu(t)$,
					$t>0$ defined above.
					The probability generating function $\mathfrak{G}^\nu (u,t) =
					\sum_{k=0}^{n_0} u^k \mathfrak{p}_k^\nu (t)$,
					$t>0$, $|u| \leq 1$, satisfies the following partial differential
					equation:
					\begin{equation}
						\label{diff}
						\frac{\partial^\nu}{\partial t^\nu}
						\mathfrak{G}^\nu (u,t)
						= \mu (n_0 +1) \left( \frac{1}{u} -1 \right)
						\left[ \mathfrak{G}^\nu (u,t) - \mathfrak{p}_0^\nu (t)
						\right] + \mu (u-1) \frac{\partial}{
						\partial u} \mathfrak{G}^\nu (u,t).
					\end{equation}
					subject to the initial condition $\mathfrak{G}^\nu (u,0) = u^{n_0}$,
					for $|u| \leq 1$, $t>0$.

					\begin{proof}
						Starting from \eqref{cna}, we obtain that
						\begin{align}
							\label{diff1}
							&\frac{d^\nu}{dt^\nu} \sum_{k=0}^{n_0}
							u^k \mathfrak{p}_k^\nu (t) \\
							& = - \mu \sum_{k=1}^{n_0} u^k
							(n_0+1-k) \mathfrak{p}_k^\nu (t)
							+ \mu \sum_{k=0}^{n_0-1} u^k (n_0-k)
							\mathfrak{p}_{k+1}^\nu (t), \notag
						\end{align}
						so that
						\begin{align}
							\frac{\partial^\nu}{\partial t^\nu}
							\mathfrak{G}^\nu (u,t) = {} & - \mu (n_0+1)
							\left[ \mathfrak{G}^\nu (u,t) - \mathfrak{p}_0^\nu (t)
							\right] + \mu u \frac{\partial}{\partial u}
							\mathfrak{G}^\nu (u,t) \\
							& + \frac{\mu (n_0+1)}{u} \left[ \mathfrak{G}^\nu (u,t)
							- \mathfrak{p}_0^\nu (t) \right] - \mu \frac{
							\partial}{\partial u} \mathfrak{G}^\nu (u,t) \notag \\
							= {} & \mu (n_0+1) \left(\frac{1}{u} -1 \right)
							\left[ \mathfrak{G}^\nu (u,t) - \mathfrak{p}_0^\nu (t)
							\right] + \mu (u-1) \frac{\partial}{\partial u}
							\mathfrak{G}^\nu (u,t). \notag
						\end{align}
					\end{proof}
				\end{theorem}

				\begin{theorem}
					The mean number of individuals $\mathbb{E} \mathfrak{M}^\nu(t)$,
					$t>0$ in the fractional sublinear death process, reads
					\begin{equation}
						\label{mean2}
						\mathbb{E} \mathfrak{M}^\nu (t) =
						\sum_{k=1}^{n_0}
						\binom{n_0+1}{k+1} (-1)^{k+1}
						E_{\nu,1} (-\mu k t^\nu), \qquad t>0, \: \nu \in (0,1].
					\end{equation}

					\begin{proof}
						From \eqref{diff} and by considering that
						$\mathbb{E} \mathfrak{M}^\nu (t)
						= \left. \frac{\partial}{\partial u}\mathfrak{G}^\nu (u,t)
						\right|_{u=1}$,
						we directly arrive at the following initial value problem:
						\begin{equation}
							\begin{cases}
								\frac{d^\nu}{dt^\nu}
								\mathbb{E} \mathfrak{M}^\nu (t) =
								- \mu (n_0+1) \left[ 1- \mathfrak{p}_0^\nu
								(t) \right] + \mu \mathbb{E} \mathfrak{M}^\nu
								(t), \\
								\mathbb{E} \mathfrak{M}^\nu (0) = n_0,
							\end{cases}
						\end{equation}
						which can be solved by resorting to the Laplace transform,
						as follows:
						\begin{align}
							\label{cons}
							\mathcal{L} \left\{ \mathbb{E}
							\mathfrak{M}^\nu (t) \right\} (z) & =
							n_0 \frac{z^{\nu-1}}{z^\nu -\mu} + \mu (n_0+1)
							\sum_{k=1}^{n_0} \binom{n_0}{k} (-1)^k
							\frac{z^{\nu-1}}{z^\nu+\mu k}
							\cdot \frac{1}{z^\nu -\mu} \\
							& = n_0 \frac{z^{\nu-1}}{z^\nu -\mu} +
							\sum_{k=1}^{n_0} \binom{n_0+1}{k+1} (-1)^k
							\left[ \frac{z^{\nu-1}}{z^\nu - \mu} -
							\frac{z^{\nu-1}}{z^\nu + \mu k} \right]. \notag
						\end{align}
						In \eqref{cons}, formula \eqref{extinction}
						must be considered.
						By inverting the Laplace transform we obtain that
						\begin{align}
							\label{inversion}
							\mathbb{E} \mathfrak{M}^\nu & (t) \\
							= {} & n_0 E_{\nu,1} (\mu t^\nu) +
							\sum_{k=1}^{n_0} \binom{n_0+1}{k+1} (-1)^k
							\left[ E_{\nu,1} (\mu t^\nu) - E_{\nu,1}
							(- \mu k t^\nu ) \right] \notag \\
							= {} & n_0 E_{\nu,1} (\mu t^\nu) +
							E_{\nu,1} (\mu t^\nu)
							\sum_{k=1}^{n_0} \binom{n_0+1}{k+1} (-1)^k \notag \\
							& - \sum_{k=1}^{n_0} \binom{n_0+1}{k+1}
							(-1)^k E_{\nu,1} (- \mu k t^\nu ) \notag \\
							{} = & \sum_{k=1}^{n_0} \binom{n_0+1}{k+1}
							(-1)^{k+1} E_{\nu,1} (- \mu k t^\nu ), \notag
						\end{align}
						as desired.
					\end{proof}
				\end{theorem}

				\begin{remark}
					The mean value \eqref{mean2} can also be directly derived as follows.
					\begin{align}
						\mathbb{E} \mathfrak{M}^\nu (t) & =
						\sum_{k=0}^{n_0} k \mathfrak{p}_k^\nu (t) \\
						& = \sum_{k=1}^{n_0} k \sum_{l=0}^{n_0-k}
						\binom{n_0-k}{l} (-1)^l E_{\nu,1} (-(l+1) \mu t^\nu ) \notag \\
						& = \sum_{k=1}^{n_0} k \sum_{l=1}^{n_0+1-k}
						\binom{n_0-k}{l-1} (-1)^{l-1} E_{\nu,1} (-\mu l t^\nu ) \notag \\
						& = \sum_{l=1}^{n_0} (-1)^{l-1} E_{\nu,1} (- \mu l t^\nu )
						\sum_{k=1}^{n_0+1-l} k \binom{n_0-k}{l-1}. \notag
					\end{align}
					It is now sufficient to show that
					\begin{equation}
						\sum_{k=1}^{n_0+1-l} k \binom{n_0-k}{l-1} =
						\binom{n_0+1}{l+1}.
					\end{equation}
					Indeed,
					\begin{align}
						\label{2mean}
						\sum_{k=1}^{n_0+1-l}
						k \binom{n_0-k}{l-1}
						& = \sum_{k=l-1}^{n_0-1} (n_0-k) \binom{k}{l-1} \\
						& = \sum_{k=l-1}^{n_0-1} (n_0+1-k-1) \binom{k}{l-1} \notag \\
						& = (n_0+1) \sum_{k=l-1}^{n_0-1} \binom{k}{l-1} -
						l \sum_{k=l-1}^{n_0-1} \binom{k+1}{l} \notag \\
						& = (n_0+1) \sum_{k=l}^{n_0} \binom{k-1}{l-1}
						- l \sum_{k=l+1}^{n_0+1} \binom{k-1}{l} \notag \\
						& = (n_0+1) \binom{n_0}{l} -l \binom{n_0+1}{l+1} \notag \\
						& = \binom{n_0+1}{l+1}. \notag
					\end{align}
					The crucial step of \eqref{2mean} is justified by the
					following formula
					\begin{equation}
						\label{formula}
						\sum_{k=j}^{n_0} \binom{k-1}{j-1} = 1 + \binom{j}{j-1}
						+ \dots + \binom{n_0-1}{j-1} = \binom{n_0}{j}.
					\end{equation}
				\end{remark}

				\begin{figure}
					\centering
					\includegraphics[scale=0.8,angle=0]{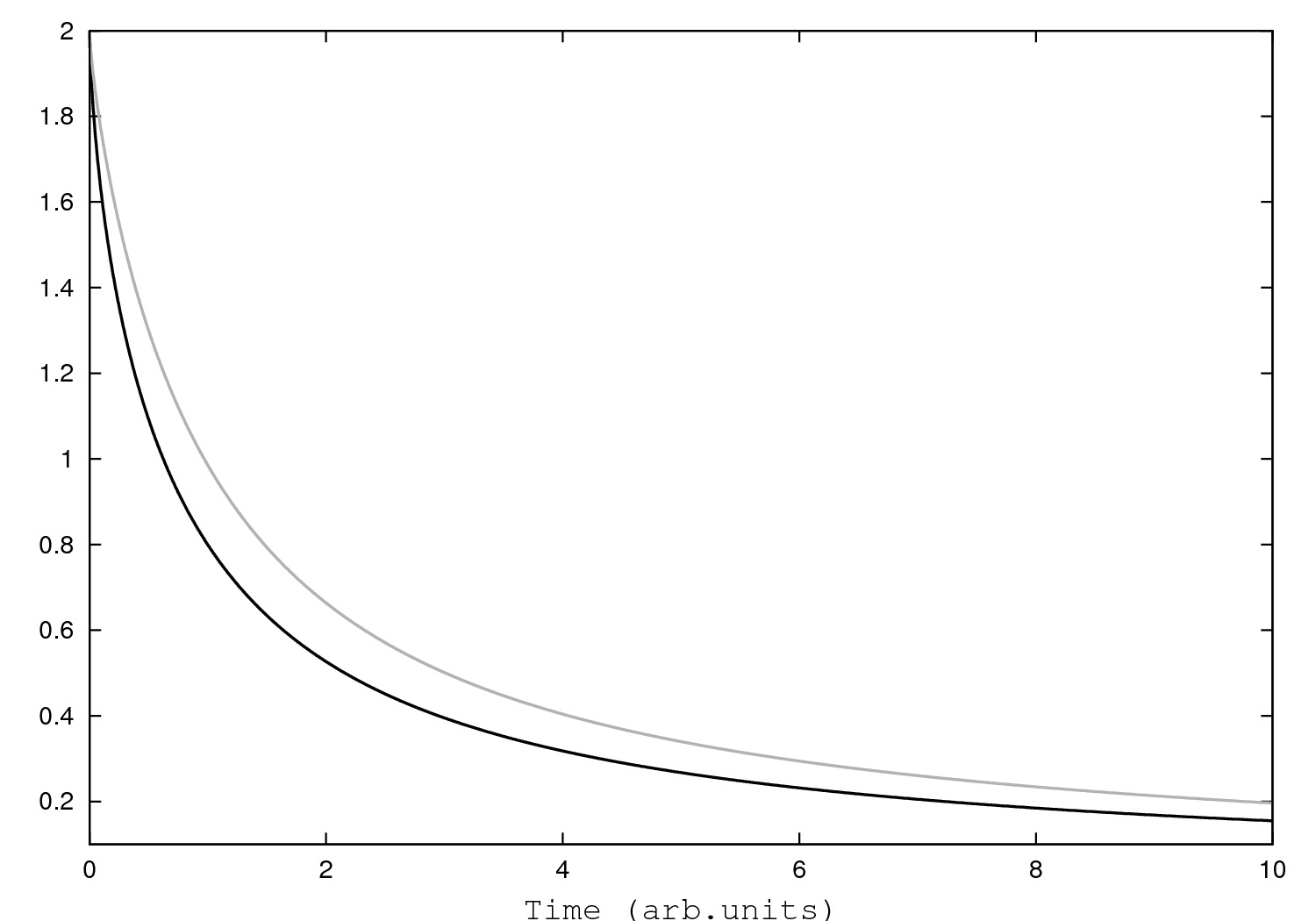}
					\caption{\label{meanplot}Plot of $\mathbb{E} M^{0.7} (t)$ (in black) and
					$\mathbb{E} \mathfrak{M}^{0.7} (t)$ (in grey), $n_0 = 2$.}
				\end{figure}

				Figure~\ref{meanplot} shows that in the sublinear case, the mean number of
				individuals in the population, decays more slowly than in the linear case,
				as expected.

				Note that \eqref{mean2} satisfies the initial condition
				$\mathbb{E} \mathfrak{M}^\nu (0) = n_0$. In order to check this, it
				is sufficient to show that
				\begin{align}
					\label{simple}
					& \sum_{k=1}^{n_0} \binom{n_0+1}{k+1} (-1)^{k+1}
					= \sum_{r=2}^{n_0+1} \binom{n_0+1}{r} (-1)^r \\
					& = \left[ \sum_{r=0}^{n_0+1} \binom{n_0+1}{r} (-1)^r \right]
					-1 + \binom{n_0+1}{1} = n_0. \notag
				\end{align}
				The details in \eqref{simple} explain also the last step of \eqref{inversion}.

			\subsubsection{Comparison of $\mathfrak{M}^\nu (t)$ with the fractional linear death
					process $M^\nu(t)$ and the fractional linear birth process
					$N^\nu (t)$}

				The distributions of the fractional linear and sublinear processes
				examined above display a behaviour which is illustrated in Table
				\ref{comparison} (see also Table 2 for the mean values).
				\begin{table}[ht]\footnotesize
					\centering
					\caption{\label{comparison}State probabilities $p_k^\nu (t)$ for the
					fractional linear death process $M^\nu (t)$, $t>0$, and
					$\mathfrak{p}_k^\nu (t)$ for the
					fractional sublinear death process $\mathfrak{M}^\nu (t)$.}
					\begin{tabular}{l}
						State Probabilities \\
						\hline \\
						$p_{n_0}^\nu (t) = E_{\nu,1} (-\mu n_0 t^\nu )$
						\vspace{0.15cm} \\
						$ \mathfrak{p}_{n_0}^\nu (t)
						= E_{\nu,1} (-\mu t^\nu)$ \\ \\
						$p_{n_0-1}^\nu (t) = n_0 \left[ E_{\nu,1}
						(- (n_0-1) \mu t^\nu)
						- E_{\nu,1} (- n_0 \mu t^\nu ) \right]$ \vspace{0.15cm} \\
						$\mathfrak{p}_{n_0-1}^\nu (t) = E_{\nu,1} (-\mu t^\nu)
						- E_{\nu,1} (-2 \mu t^\nu)$ \\
						\vdots \\
						$p_k^\nu (t) = \binom{n_0}{k} \sum_{l=0}^{n_0-k}
						\binom{n_0-k}{l} (-1)^l E_{\nu,1} (- (k+l) \mu t^\nu )$
						\vspace{0.15cm} \\
						$\mathfrak{p}_k^\nu (t) = \sum_{l=0}^{n_0-k}
						\binom{n_0-k}{l} (-1)^l	E_{\nu,1}
						\left( - \left( l+1 \right)
						\mu t^\nu \right)$ \\
						\vdots \\
						$p_1^\nu (t) = n_0 \sum_{l=0}^{n_0-1}
						\binom{n_0-1}{l} (-1)^l E_{\nu,1} (- (1+l) \mu t^\nu )$
						\vspace{0.15cm} \\
						$\mathfrak{p}_1^\nu (t)
						= \sum_{l=0}^{n_0-1}
						\binom{n_0-1}{l} (-1)^l	E_{\nu,1}
						\left( - \left( l+1 \right)
						\mu t^\nu \right)$ \\ \\
						$p_0^\nu (t) = \sum_{l=0}^{n_0}
						\binom{n_0}{l} (-1)^l E_{\nu,1} (-l \mu t^\nu )$
						\vspace{0.15cm} \\
						$\mathfrak{p}_0^\nu (t) =
						\sum_{l=0}^{n_0} \binom{n_0}{l} (-1)^l E_{\nu,1}
						(-l \mu t^\nu)$ \\ \\
						\hline
					\end{tabular}
				\end{table}

				The most striking fact about the models dealt with above, is that the
				linear probabilities decay faster than the corresponding
				sublinear ones, for small values of $k$; whereas, for
				large values of $k$, the sublinear probabilities
				take over and the extinction probabilities in both cases coincide.
				The reader should also compare the state probabilities of the death models
				examined here with those of the fractional linear pure birth process
				(with birth rate $\lambda$ and one progenitor). These read
				\begin{equation}
					\hat{p}^\nu_k (t) = \sum_{j=1}^k \binom{k-1}{j-1}
					(-1)^{j-1} E_{\nu,1} (- \lambda j t^\nu), \qquad k \geq 1.
				\end{equation}
				Note that $\hat{p}_1^\nu (t) = E_{\nu,1} (-\lambda t^\nu)$ is of the
				same form as $\mathfrak{p}_{n_0}^\nu (t) = E_{\nu,1} (-\mu t^\nu)$.
				We now show that
				\begin{align}
					\label{form}
					\sum_{k=n_0+1}^\infty \hat{p}_k^\nu (t) & =
					1 - \sum_{k=1}^{n_0} \hat{p}_k^\nu (t) \\
					& = 1 - \sum_{k=1}^{n_0} \sum_{j=1}^k \binom{k-1}{j-1}
					(-1)^{j-1} E_{\nu,1} (- \lambda j t^\nu) \notag \\
					& = 1 - \sum_{j=1}^{n_0} (-1)^{j-1}
					E_{\nu,1} (- \lambda j t^\nu) \sum_{k=j}^{n_0} \binom{k-1}{j-1}
					\notag \\
					& = 1 - \sum_{j=1}^{n_0} (-1)^{j-1} \binom{n_0}{j}
					E_{\nu,1} (- \lambda j t^\nu) \notag \\
					& = \text{\eqref{extinction} with $\lambda$ replacing $\mu$}. \notag
				\end{align}
				Note that in the above step we used formula \eqref{formula}.

				\begin{table}\footnotesize
					\centering
					\caption{\label{mean-table}Mean values for the
						fractional linear birth $N^\nu (t)$,
						fractional linear death $M^\nu (t)$ and fractional sublinear
						death $\mathfrak{M}^\nu (t)$ processes.}
					\begin{tabular}{l}
						\hline \\
						$ \mathbb{E}N^\nu (t) = E_{\nu,1} (\lambda t^\nu)$ \\ \\
						$ \mathbb{E}M^\nu (t) = n_0 E_{\nu,1} (- \mu t^\nu)$ \\ \\
						$ \mathbb{E}\mathfrak{M}^\nu (t) = \sum_{k=1}^{n_0}
						\binom{n_0+1}{k+1} (-1)^{k+1} E_{\nu,1} (-\mu k t^\nu)$ \\ \\
						\hline \\
					\end{tabular}
				\end{table}

				By comparing formulae (3.4) of \ocite{pol} and \eqref{pik} above, we
				arrive at the conclusion that (for $\lambda=\mu$)
				\begin{align}
					& \text{Pr} \left\{ N^\nu (t) = k \mid N^\nu (0) = 1 \right\} \\
					& = \sum_{j=1}^k \binom{k-1}{j-1} (-1)^{j-1} E_{\nu,1}
					(- \lambda_j t^\nu ) \notag \\
					& = \text{Pr} \left\{ \mathfrak{M}^\nu (t) = n_0 +1 - k
					\mid \mathfrak{M} (0) = n_0 \right\}, \qquad 1 \leq k \leq n_0.
					\notag
				\end{align}
				The probability of extinction $\mathfrak{p}_0^\nu(t)$ corresponds to
				the probability of the event $\left\{ N^\nu (t) > n_0 \right\}$
				for the fractional linear birth process.

				\vspace{.2cm}
				\textit{Acknowledgement:} The authors wish to thank Francis Farrelly for having
					checked and corrected the manuscript. The authors are grateful to
					the referees for drawing our attention to some relevant references
					and for useful remarks which improved the presentation of the paper.

\end{document}